\documentclass[11pt,reqno]{amsart}    
\usepackage{amsmath,amsfonts,amsthm,amsopn,amssymb}
\usepackage{latexsym,soul,cite,mathrsfs}
\usepackage{color,xcolor}
\usepackage{color,enumitem,graphicx}
\usepackage[colorlinks=true,urlcolor=blue,
citecolor=red,linkcolor=blue,linktocpage,pdfpagelabels,
bookmarksnumbered,bookmarksopen]{hyperref}
\usepackage[english]{babel}
\usepackage{enumitem}
\usepackage[left=3.0cm,right=3.0cm,top=2.8cm,bottom=2.8cm]{geometry}

\usepackage[hyperpageref]{backref}
\makeatletter
\newcommand{\leqnomode}{\tagsleft@true}
\newcommand{\reqnomode}{\tagsleft@false}
\makeatother

\numberwithin{equation}{section}

\newtheorem{thm}{Theorem}[section]
\newtheorem{lem}[thm]{Lemma}

\title[Compactness of the Extremal Functions]{Compactness of Extremals for Singular Anisotropic Trudinger-Moser functionals on bounded domain}

\author{Weiwei Shan, Minbo Yang,  Jiazheng Zhou}

\address{Weiwei Shan  \newline\indent School of Mathematical Sciences, Zhejiang Normal University, \newline\indent
	Jinhua, Zhejiang, 321004, People's Republic of China}
\email{weiweishan@zjnu.edu.cn}

\address{Minbo Yang  \newline\indent School of Mathematical Sciences, Zhejiang Normal University, \newline\indent
Jinhua, Zhejiang, 321004, People's Republic of China}
\email{mbyang@zjnu.edu.cn}

\address{Jiazheng Zhou  \newline\indent  Departamento de Matem\'{a}tica, Universidade de Bras\'{i}lia\newline\indent Bras\'{i}lia, DF 70910-900, Brasil}
\email{zhou@mat.unb.br} 


 \thanks{}


\begin{document}
	\maketitle
\begin{abstract}
In this paper, we investigate the compactness of extremal functions
for a critical singular anisotropic Trudinger-Moser inequality established by Lu-Shen-Xue-Zhu\cite{ref1}. We prove by means of blow-up analysis that the extremals $u_{\beta}$ converge  in $W_{0}^{1,n}(\Omega)\cap C^{1}(\overline{\Omega})$ to some function $u_{0}$ which achieves the supremum
 \begin{equation}
 \sup\limits_{u\in W_{0}^{1,n}(\Omega),\Vert u\Vert_{F(\Omega)}\leq1}\int_{\Omega}^{}e^{\tau_{n}\vert u\vert^{\frac{n}{n-1}}}dx,\notag
 \end{equation}
as $\beta\to 0$, where $\tau_{n}=n^{\frac{n}{n-1}}\kappa_{n}^{\frac{1}{n-1}}$, $\kappa_{n}$ denotes the volume of the unit Wulff ball in $\mathbb{R}^{n}$ and $\Vert u\Vert_{F(\Omega)}$ is the anisotropic norm of $u$.
\end{abstract}

\section{Introduction And Main Results}
Let $\Omega\subset\mathbb{R}^{n}(n\geq2)$ be a smooth bounded domain, $ W_{0}^{1,n}(\Omega) $ be the usual Sobolev space, that is, the completion $C_{0}^{\infty}(\Omega)$ equipped with the norm $\Vert u\Vert_{W_{0}^{1,n}(\Omega)}=(\int_{\Omega}^{}\vert \nabla u\vert^{n}dx)^{1/n}$. Then Sobolev embedding theorem states that embedding $W_{0}^{1,q}(\Omega)\hookrightarrow L^{\frac{nq}{n-q}}(\Omega)$ holds for all $q\in\left[1,n\right)$. But when $q=n$, it
is known that $W_{0}^{1,n}(\Omega)\not\hookrightarrow L^{\infty}(\Omega)$. In this case, Yudovich \cite{ref2}, Pohozaev \cite{ref3}, Peetre \cite{ref4} and Trudinger \cite{ref5} independently proved that the embedding
$W_{0}^{1,n}(\Omega)\hookrightarrow L^{\phi}(\Omega)$ is optimal, where $L^{\phi}(\Omega)$ denotes the Orlicz space associated with the function
$
\phi(t)=e^{\alpha\vert t\vert^{\frac{n}{n-1}}}-1$
for some $\alpha>0$. Later, Moser\cite{ref6} employed the rearrangement theory to find the optimal exponent $\alpha$.
The classical Trudinger-Moser inequality states the following:
\begin{equation}\label{E:1.1}
\sup\limits_{u\in W_{0}^{1,n}(\Omega),\Vert u \Vert_{W_{0}^{1,n}(\Omega)}\leq1}\int_{\Omega}^{}e^{\alpha \vert u\vert^{\frac{n}{n-1}}}dx\left\{  
\begin{aligned}
&<+\infty,\quad\forall\alpha\leq \alpha_{n},\\
&=+\infty,\quad\forall\alpha>\alpha_{n},\\
	\end{aligned}
	\right.
\end{equation}
where $\alpha_{n}=n^{\frac{n}{n-1}}V_{n}^{\frac{1}{n-1}}$, $V_{n}=\vert\mathbb{B} \vert$ and $\mathbb{B}:=\left\{x\in\mathbb{R}^{n}:\vert x\vert\leq1\right \} $. Notably, for $\alpha\leq \alpha_{n}$, the supremum in \eqref{E:1.1} can be attained by some function. This result was first established by Carleson-Chang\cite{ref7} for the case $\Omega=\mathbb{B}$.  Precisely, there exists some function 
$u_{\beta}\in W_{0}^{1,n}(\mathbb{B})$ with $\Vert u_{\beta} \Vert_{W_{0}^{1,n}(\mathbb{B})}=1$ such that
\begin{equation}
\int_{\mathbb{B}}^{}e^{\alpha \vert u_{\beta}\vert^{\frac{n}{n-1}}}dx=\sup\limits_{u\in W_{0}^{1,n}(\mathbb{B}),\Vert u \Vert_{W_{0}^{1,n}(\mathbb{B})}\leq1}\int_{\mathbb{B}}^{}e^{\alpha \vert u\vert^{\frac{n}{n-1}}}dx.\notag
\end{equation} 
Struwe\cite{ref8} subsequently showed that when $\Omega$ and $\mathbb{B}$ are close in measure, the supremum in \eqref{E:1.1} can be achieved. For any bounded domain $\Omega\in\mathbb{R}^{2}$, Flucher\cite{ref9} obtained similar results using conformal mapping techniques. Then Lin\cite{ref10} extended Flucher's results to higher-dimensional spaces.

Inequality \eqref{E:1.1} has been extended in various forms. By a rearrangement approach,  Adimurthi-Sandeep \cite{ref11} established the singular form of inequality \eqref{E:1.1}: for any $\beta\in(0,1)$, there holds
\begin{equation}\label{E:1.2}
\sup\limits_{u\in W_0^{1,n}(\Omega),\Vert u \Vert_{W_{0}^{1,n}(\Omega)}\leq1}\int_\Omega \frac{e^{\alpha(1-\beta)\vert u\vert^{\frac{n}{n-1}}}}{|x|^{n\beta}}dx\left\{  
\begin{aligned}
&<+\infty,\quad\forall\alpha\leq \alpha_{n},\\
&=+\infty,\quad\forall\alpha> \alpha_{n}.\\
	\end{aligned}
	\right.
\end{equation}
Csat\'{o}-Roy\cite{ref12}
and Yang-Zhu\cite{ref13} independently proved the existence of extremal functions for the above inequality in the case $n=2$, while the higher-dimensional cases $(n\geq3)$ were proved by Csat\'{o}-Roy-Nguyen\cite{ref14}.

Since the embedding $W_{0}^{1,n}(\mathbb{B})\cap\ell\hookrightarrow L^{p}(\mathbb{B},\vert x\vert^{k})$ holds for all $k>0$ and $p=2(n+k)/(n-2)$, where $\ell$ denotes the set of all radial functions. Based on the work of de Figueiredo-do \'{O}-dos Santos\cite{ref42}, Yang-Zhu\cite{ref43} established the weighted form of inequality \eqref{E:1.1} as follows:
\begin{equation}\label{E:1.333}
	\sup\limits_{u\in W_0^{1,n}(\mathbb{B})\cap\ell,\Vert u \Vert_{W_{0}^{1,n}(\mathbb{B})}\leq1}\int_{\mathbb{B}}|x|^{n\varepsilon}e^{\alpha(1+\varepsilon)\vert u\vert^{\frac{n}{n-1}}}dx\left\{  
	\begin{aligned}
		&<+\infty,\quad\forall\alpha\leq \alpha_{n},\\
		&=+\infty,\quad\forall\alpha> \alpha_{n},\\
	\end{aligned}
	\right.
\end{equation}
for any $\varepsilon>0$. Moreover, the above supremum is attained by some function $u_{\varepsilon}\in W_0^{1,n}(\mathbb{B})\cap\ell$ with $\Vert u \Vert_{W_{0}^{1,n}(\mathbb{B})}=1$.

A natural problem is to investigate the compactness of extremal functions for the critical Trudinger-Moser inequality. The first result on compactness was established by Wang-Yang\cite{ref15}, they proved the following: Let $\{u_{\beta}\}$ be a sequence of extremal functions for the critical inequality \eqref{E:1.2} in the case $n=2$, where each $u_{\beta}\in W_{0}^{1,2}(\Omega)$ satisfies $\Vert u_{\beta}\Vert_{W_{0}^{1,2}(\Omega)}=1$ and 
\begin{equation}
	\int_\Omega \frac{e^{4\pi(1-\beta) u_{\beta}^{2}}}{|x|^{2\beta}}dx=\sup\limits_{u\in W_0^{1,2}(\Omega),\Vert u \Vert_{W_{0}^{1,2}(\Omega)}\leq1}\int_\Omega \frac{e^{4\pi(1-\beta) u^{2}}}{|x|^{2\beta}}dx.\notag
\end{equation}  Then by passing to a subsequence, $u_{\beta}\to u_{0}$ in $C^{1}(\overline{\Omega})$ as $\beta\to0$. Moreover, the limit function 
$u_{0}$ achieves the supremum 
\begin{equation}
\sup\limits_{u\in W_0^{1,2}(\Omega),\Vert u \Vert_{W_{0}^{1,2}(\Omega)}\leq1}\int_\Omega e^{4\pi u^{2}}dx.\notag
\end{equation}
Su-Xie-Li\cite{ref16} extended the above results to higher-dimensional cases. For the weighted inequality \eqref{E:1.333} in the case $n=2$, Shan-Li\cite{ref17} recently proved that, up to a subsequence, $u_{\varepsilon}\to u_{0}$ in $C^{1}(\overline{\mathbb{B}})$ as $\varepsilon\to0$ and $u_{0}$ is also an extremal function of the classical Trudinger-Moser inequality. A similar result holds for the singular Hardy-Trudinger-Moser inequality(see\cite{ref18}). For results on unbounded domains, see Li\cite{ref19}.
 
In the present paper we are interested in 
the anisotropic Trudinger-Moser inequality.
Let $F(x)\in C^{2}(\mathbb{R}^{n}\setminus\left \{0\right\})$ be a nonnegative even function in $\mathbb{R}^{n}$. Furthermore, $F(x)$ is a convex function satisfying
$$
F(a\xi)=\vert a\vert F(\xi),~\forall~\xi\in \mathbb{R}^{n},a\in\mathbb{R}.\notag
$$
Since $F(x)$ is a homogeneous function, there exist two constants $0<a_{1}\leq a_{2}<+\infty$ satisfying
$a_{1}\vert \xi\vert\leq F(\xi)\leq a_{2}\vert \xi\vert.$ A typical example satisfying the above conditions is $$F(\xi)=({\textstyle \sum_{i}}\vert \xi_{i}\vert^{p})^{\frac{1}{p}},~\forall~p\in\left[ 1,+\infty\right ).$$
Clearly, when $p=2$, 
$F(x)$ reduces to the Euclidean norm 
$\Vert x\Vert$. Notably, Hess($F^{n}$) is positive definite in $\mathbb{R}^{n}\setminus\{0\}$(see\cite{ref44}).
The anisotropic Dirichlet norm of $u$ is defined by 
\begin{equation}
\Vert u \Vert_{F(\Omega)}=\bigg(\int_{\Omega}^{}F^{n}(\nabla u)dx\bigg)^{1/n}\notag
\end{equation}
with Euler equation involves the operator
\begin{equation}
\mathrm{Q}_{n}u:= {\textstyle \sum_{i=1}^{i=n}}\frac{\partial }{\partial x_{i}}(F^{n-1}(\nabla u)F_{\xi_{i}}(\nabla u)).\notag
\end{equation}
When $F(x)\not=\Vert x\Vert$, this operator becomes nonlinear and is called the n-Finsler-Laplacian.

Let $\mathbb{S}^{n-1}$ be the boundary of $\mathbb{B}$. The map $\mathcal{M}:\mathbb{S}^{n-1}\mapsto\mathbb{R}^{n},$ defined by $\mathcal{M}(\xi)=F_{\xi}(\xi)$ where $F_{\xi}=(F_{\xi_{1}},F_{\xi_{2}},\dots,F_{\xi_{n}})$. Its image $\mathcal{M}(\mathbb{S}^{n-1})$ is a smooth,
convex hypersurface in $\mathbb{R}^{n}$, which is called Wulff shape of $F$. The support function of the set $\mathbb{A}=\left\{x\in\mathbb{R}^{n}:F(x)\leq1\right\}$ is given by $F^{\circ}(x)=\sup\limits_{\xi\in K}\left\langle x,\xi\right\rangle$. Furthermore, $F^{\circ}(x)\in C^{2}(\mathbb{R}^{n}\setminus\left \{0\right\})$ is a convex homogeneous function and the dual function of $F(x)$, satisfying the duality relations:
$$F(x)=\sup\limits_{\xi\in\mathbb{R}^{n}\setminus\left \{0\right\}}\frac{\left\langle  x,\xi\right\rangle}{F^{\circ}(\xi)},\quad\quad F^{\circ}(x)=\sup\limits_{\xi\in\mathbb{R}^{n}\setminus\left\{0\right\}}\frac{\left\langle x,\xi\right\rangle}{F(\xi)}. $$
Clearly, $\mathcal{M}(\mathbb{S}^{n-1})=\left\{x\in\mathbb{R}^{n}:F^{\circ}(x)=1\right\}.$ 
Define the Wulff ball centered at $x_{1}$ with radius $R>0$ by $\mathcal{W}_{R}(x_{1}):=\left\{x\in\mathbb{R}^{n}:F^{\circ}(x-x_{1})\leq R\right\}.$

By replacing $\Vert u \Vert_{W_{0}^{1,n}(\Omega)}$ with $\Vert u \Vert_{F(\Omega)}$ in \eqref{E:1.1}, Wang-Xia \cite{ref20} applied convex rearrangement theory(see\cite{ref21}) to establish an anisotropic
Trudinger-Moser inequality: 
\begin{equation}\label{E:1.3}
\sup\limits_{u\in W_{0}^{1,n}(\Omega),\Vert u \Vert_{F(\Omega)}\leq1}\int_{\Omega}^{}e^{\tau \vert u\vert^{\frac{n}{n-1}}}dx\left\{  
\begin{aligned}
&<+\infty,\quad\forall\tau\leq \tau_{n},\\
&=+\infty,\quad\forall\tau>\tau_{n},\\
	\end{aligned}
	\right.
\end{equation}
where $\tau_{n}=n^{\frac{n}{n-1}}\kappa_{n}^{\frac{1}{n-1}}$ and $\kappa_{n}=\vert\mathcal{W}_{1}\vert$. By using blow-up analysis, Zhou-Zhou \cite{ref22} showed that the supremum in \eqref{E:1.3} is attainable for $\tau\leq \tau_{n}$.
Similar to the Adimurthi-Sandeep's result \eqref{E:1.2},  Zhu\cite{ref23} obtained 
the following inequality:
\begin{equation}\label{E:1.4}
\sup\limits_{u\in W_{0}^{1,n}(\Omega),\Vert u \Vert_{F(\Omega)}\leq1}\int_{\Omega}^{} \frac{e^{\tau(1-\beta)\vert u\vert^{\frac{n}{n-1}}}}{F^{\circ}(x)^{n\beta}}dx\left\{ \begin{aligned}
&<+\infty,\quad\forall\tau\leq \tau_{n},\\
&=+\infty,\quad\forall\tau> \tau_{n},\\
	\end{aligned}
	\right.
\end{equation}
for any $\beta\in(0,1)$. Recently, Lu-Shen-Xue-Zhu \cite{ref1} proved the existence of extremal functions for \eqref{E:1.4} when $\tau\leq \tau_{n}$. Specifically, in the critical case $\tau=\tau_{n}$, there exists some function $u_{\beta}\in W_{0}^{1,n}(\Omega)$ with $\Vert u_{\beta}\Vert_{F(\Omega)}=1$ satisfying
\begin{equation}\label{E:1.5}
\int_{\Omega}^{} \frac{e^{\tau_{n}(1-\beta)\vert u_{\beta}\vert^{\frac{n}{n-1}}}}{F^{\circ}(x)^{n\beta}}dx=\sup\limits_{u\in W_{0}^{1,n}(\Omega),\Vert u\Vert_{F(\Omega)}\leq1}\int_{\Omega}^{} \frac{e^{\tau_{n}(1-\beta)\vert u\vert^{\frac{n}{n-1}}}}{F^{\circ}(x)^{n\beta}}dx.
\end{equation}
For other results related to the anisotropic Trudinger-Moser inequality, we refer to \cite{ref24,ref25,ref26,ref27,ref28}.

Inspired by the work of Wang-Yang\cite{ref15}, we aim to investigate the compactness of the  extremals for critical inequality \eqref{E:1.4} in this paper. For simplicity, let 
\begin{equation}\label{FBS}
\mathfrak{T}=\sup\limits_{u\in W_{0}^{1,n}(\Omega),\Vert u\Vert_{F(\Omega)}\leq1}\int_{\Omega}^{} e^{\tau_{n}\vert u\vert^{\frac{n}{n-1}}}dx.
\end{equation}
We now present the main result.
\begin{thm}\label{th:1.1}
Let $0\in\Omega$ and $\left \{u_{\beta}\right \} $ be the extremal function sequence satisfying $\Vert u_{\beta}\Vert_{F(\Omega)}=1$ and
\eqref{E:1.5} for $0<\beta<1$. Then as $\beta\to 0$, passing to a subsequence, $u_{\beta}\to u_{0}$ in $ W_{0}^{1,n}(\Omega)\cap C^{1}(\overline{\Omega})$. Moreover, $u_{0}$ satisfies \begin{equation}
\int_{\Omega}^{}e^{\tau_{n}\vert u_{0}\vert^{\frac{n}{n-1}}}dx=\mathfrak{T}.\notag
\end{equation}
\end{thm}

The proof of Theorem \ref{th:1.1} relies crucially on the blow-up analysis that was originally developed by Ding-Jost-Li-Wang \cite{ref29} and Li \cite{ref30,ref31}. And it has been widely applied in establishing the existence of extremal functions for Trudinger-Moser type inequalities (see \cite{ref32,ref33,ref34,ref35}). Suppose $u_{\beta}$ is given as in \eqref{E:1.5} and it is also the solution of equation \eqref{E:3.1} below. If blow-up occurs for $u_{\beta}$ in $\Omega$, we can obtain $\mathfrak{T}\leq\vert\Omega\vert+\kappa_{n}e^{\tau_{n}C_{x_{0}}+{\sum_{j=1}^{n-1}}\frac{1}{j}}$. However,  we can find a test function $\upsilon_{\beta}(x)\in W_{0}^{1,n}(\Omega)$ satisfying $\Vert\upsilon_{\beta}(x)\Vert_{F(\Omega)}=1$ and 
\begin{equation*}
	\int_{\Omega}^{} e^{\tau_{n}\vert\upsilon_{\beta}(x)\vert^{\frac{n}{n-1}}}dx>\vert\Omega\vert+\kappa_{n}e^{\tau_{n}C_{x_{0}}+{\sum_{j=1}^{n-1}}\frac{1}{j}},
\end{equation*}
which implies 
\begin{equation}
	\mathfrak{T}>\vert\Omega\vert+\kappa_{n}e^{\tau_{n}C_{x_{0}}+{\sum_{j=1}^{n-1}}\frac{1}{j}}.\notag
\end{equation} Therefore, we can prove that the blow-up doesn't happen and so $\Vert u_{\beta}\Vert_{L^{\infty}(\Omega)}<+\infty$. Notably, the problem involves the n-Finsler-Laplacian, which needs to be combined with the Wulff shape of $F(x)$ in the subsequent analysis.

The paper is organized as follows: In Section 2, we first introduce several rearrangement methods and review relevant existing results. Subsequently, we apply variational methods to deduce the Euler-Lagrange equation for 
$u_{\beta}$. In Section 3, we establish an upper bound estimate for 
$\mathfrak{T}$ by blow-up analysis. In Section 4, we eliminate the blow-up
phenomenon by constructing a test function, which proves the theorem.

\section{Notations and Preliminaries}

In what follows, we denote constants by $C$. The notation $ o_{\beta}(1) $ represents an infinitesimal depending on $\beta$ such that $ o_{\beta}(1)\to 0$ as $\beta\to0$. Similarly, $ o_{\beta,R}(1)\rightarrow0 $ as $ \beta\rightarrow0 $ for any fixed $ R>0 $.

Let $u$ be a measurable function  on $\Omega\in\mathbb{R}^{n}$. Define
\begin{equation}
	u^{\ast}=\sup\left\{s\geq0:\vert\left\{x\in\Omega:\vert u(x)\vert>s\right\}\vert>t\right\},~\mathrm {for}~ t\in\mathbb{R},\notag
\end{equation}
as the one-dimensional decreasing rearrangement of $u$. Assume $\vert\mathbb{B}_{R_{1}}\vert=\vert\Omega\vert=\vert\mathcal{W}_{R_{2}}\vert$, we define 
\begin{equation}\label{E:2.2}
	u^{\sharp}=u^{\ast}(V_{n}\vert x\vert^{n}),~\mathrm {for}~ x\in\mathbb{B}_{R_{1}}\mathrm {(Schwarz~ rearrangement)}
\end{equation}
and
\begin{equation}\label{E:2.1}
	u^{\star}=u^{\ast}(\kappa_{n}F^{\circ}(x)^{n}),~\mathrm {for}~ x\in\mathcal{W}_{R_{2}}\mathrm {(Convex~ rearrangement)}.
\end{equation}
Clearly, the convex symmetrization(see\cite{ref21}) generalizes the schwarz symmetrization(see\cite{ref51}).

We now outline some properties of the functions $F(x)$ and $F^{\circ}(x)$(see\cite{ref36,ref37}): 
\begin{lem}\label{le:2.1}There holds\\
(1).$\left\langle x,\nabla F^{\circ}(x) \right\rangle=F^{\circ}(x),\left\langle x,\nabla F(x) \right\rangle=F(x)$ for any $x\neq0$;\\
(2).$\vert F(y)-F(x)\vert\leq F(y+x)\leq F(y)+F(x)$;\\
(3).$\frac{1}{C}\leq\vert F^{\circ}(x)\vert\leq C,\frac{1}{C}\leq\vert F(x)\vert\leq C$ for some $C>0$ and $x\neq0$;\\
(4).$F^{\circ}(\nabla F(x))=1,F(\nabla F^{\circ}(x))=1$ for $x\neq0$;\\
(5).$F_{\xi}(t\xi)=sgn(t)F_{\xi}(\xi)$ for $t\neq0$ and $x\neq0$;\\
(6).$F(x)F_{\xi}^{\circ}(\nabla F(x))=x=F^{\circ}(x)F_{\xi}(\nabla F^{\circ}(x))$ for $x\neq0$.
\end{lem}

From Theorem 1.2 of \cite{ref27}, we have the following Lions type concentration-compactness principle for the singular anisotropic Trudinger-Moser inequality: 
\begin{lem}\label{le:2.2}
Let $\left\{ u_{k}\right\}$ is a sequence in $ W_{0}^{1,n}(\Omega)$ satisfying $\Vert u_{k}\Vert_{F(\Omega)}=1$. If $u_{k}\rightharpoonup u_{0}\not\equiv0$ weakly in $W_{0}^{1,n}(\Omega)$, then for $0<\beta<1$,  $0<p<1/(1-\Vert u_{0}\Vert_{F(\Omega)}^{n})^{1/(n-1)}$, there holds
	$$\int_{\Omega}^{} \frac{e^{\tau_{n}(1-\beta)p\vert u_{k}\vert^{\frac{n}{n-1}}}}{F^{\circ}(x)^{n\beta}}dx<+\infty.$$ 
	This assertion becomes invalid when $p\geq1/(1-\Vert u_{0}\Vert_{F(\Omega)}^{n})^{1/(n-1)}$.
\end{lem}

According to Theorems 2.2 and 4.5 in \cite{ref22}, we have the following two estimates:
\begin{lem}\label{le:2.3}
Suppose $u\in W_{0}^{1,n}(\Omega)$ is a solution of
	\begin{equation}
		-\mathrm{Q}_{n}u=f.\notag
	\end{equation}
For any $p>1$, the estimate $$\Vert u\Vert_{L^{\infty}(\Omega)}\leq C\Vert f\Vert_{L^{p}(\Omega)}^{\frac{1}{n-1}}$$ holds whenever $f\in L^{p}(\Omega)$, where C is a fixed constant.
\end{lem}

\begin{lem}\label{le:2.4}
Assume that $\{f_{\beta}\}\in L^{1}(\Omega)$, and  $\{u_{\beta}\}\subset W_{0}^{1,n}(\Omega)\cap C^{1}(\overline{\Omega})$ satisfies
\begin{equation}
-\mathrm{Q}_{n}u_{\beta}=f_{\beta} \quad\quad\mathrm{in}\quad\Omega.\notag
\end{equation}
Then for $1<p<n$, the estimate $$\Vert \nabla u_{\beta}\Vert_{L^{p}(\Omega)}\leq C\Vert f_{\beta}\Vert_{L^{1}(\Omega)}$$ holds, where $C$ is a fixed constant.
\end{lem}

For all $\beta\in(0,1)$, we set
\begin{equation}
\mathfrak{T}_{\beta}=\sup\limits_{u\in W_{0}^{1,n}(\Omega),\Vert u\Vert_{F(\Omega)}\leq1}\int_{\Omega}^{} \frac{e^{\tau_{n}(1-\beta)\vert u\vert^{\frac{n}{n-1}}}}{F^{\circ}(x)^{n\beta}}dx=\int_{\Omega}^{} \frac{e^{\tau_{n}(1-\beta)\vert u_{\beta}\vert^{\frac{n}{n-1}}}}{F^{\circ}(x)^{n\beta}}dx,\notag
\end{equation}
where $ u_{\beta}\in W_{0}^{1,n}(\Omega)$ and $\Vert u_{\beta}\Vert_{F(\Omega)}=1 $.  By the existence result in \cite{ref1}, We know that there exists the extremal function $ u_{\beta} $ solves 
\begin{equation}\label{E:3.1}
-\mathrm{Q}_{n} u_{\beta}=\frac{1}{\lambda_{\beta}}F^{\circ}(x)^{-n\beta}u_{\beta}^{\frac{1}{n-1}}e^{\tau_{n}(1-\beta)u_{\beta}^{\frac{n}{n-1}}}~\mathrm{in}~\Omega, 
\end{equation}
where $u_{\beta}\geq0$ if and only if $u_{\beta}=0$ on $\partial\Omega$, $\Vert u_{\beta}\Vert_{F(\Omega)}=1 $, and 
\begin{equation}
\lambda_{\beta}=\int_{\Omega}^{}F^{\circ}(x)^{-n\beta}u_{\beta}^{\frac{n}{n-1}}e^{\tau_{n}(1-\beta)u_{\beta}^{\frac{n}{n-1}}}dx.\notag
\end{equation}
The Lagrange multiplier $\lambda_{\beta}$ is a positive constant.  In fact,
 the boundedness of $u_{\beta} $ in $W_{0}^{1,n}(\Omega) $ implies the existence of some function $u_{0}$ such that
\begin{align}
&u_{\beta}\rightharpoonup u_{0} \quad \mathrm {weakly ~in}~ W_{0}^{1,n}(\Omega),\notag\\
&u_{\beta}\rightarrow u_{0} \quad \mathrm {strongly ~in}~L^{s}(\Omega),~\forall~s>1 ,\notag\\
&u_{\beta}\rightarrow u_{0} \quad\mathrm {a.e.~in}~\Omega.\notag
\end{align}
\begin{lem}\label{le:3.1}
There holds
\begin{equation}
\liminf\limits_{\beta\rightarrow 0 }\lambda_{\beta}>0.\notag
\end{equation}
\end{lem}
\begin{proof}
For any $x\in \mathcal{W}_{1}$ and $\beta\in\left [0,b\right )$, where $b>0$ is sufficiently small, we have
\begin{equation}
0<\frac{e^{\tau_{n}(1-\beta)\vert u\vert^{\frac{n}{n-1}}}}{F^{\circ}(x)^{n\beta}}\leq \frac{e^{\tau_{n}\vert u\vert^{\frac{n}{n-1}}}}{F^{\circ}(x)^{bn}}.\notag
\end{equation}
Obviously, $e^{\tau_{n}\vert u\vert^{\frac{n}{n-1}}}/F^{\circ}(x)^{bn}\in L^{1}(\mathcal{W}_{1})$ for fixed $u\in W_{0}^{1,n}(\mathcal{W}_{1})$. On the flip side, for any $x\in \Omega\setminus \mathcal{W}_{1}$, one sees that
\begin{equation}
0<\frac{e^{\tau_{n}(1-\beta)\vert u\vert^{\frac{n}{n-1}}}}{F^{\circ}(x)^{n\beta}}\leq e^{\tau_{n}\vert u\vert^{\frac{n}{n-1}}}.\notag
\end{equation}
It follows from \eqref{E:1.3} that $e^{\tau_{n}\vert u\vert^{\frac{n}{n-1}}}\in L^{1}(\Omega\setminus \mathcal{W}_{1})$.
The Lebesgue Dominated Convergence Theorem leads to
\begin{equation}\label{E:3.2}
	\int_{\Omega}^{} e^{\tau_{n}\vert u\vert^{\frac{n}{n-1}}}dx= \lim\limits_{\beta\rightarrow0}\int_{\Omega}^{} \frac{e^{\tau_{n}(1-\beta)\vert u\vert^{\frac{n}{n-1}}}}{F^{\circ}(x)^{n\beta}}dx.\notag
\end{equation}
It is easy to see that
\begin{align}
\int_{\Omega}^{} e^{\tau_{n}\vert u\vert^{\frac{n}{n-1}}}dx
\leq\liminf\limits_{\beta\rightarrow0}\mathfrak{T}_{\beta}=\liminf\limits_{\beta\rightarrow0}\int_{\Omega}^{} \frac{e^{\tau_{n}(1-\beta) u_{\beta}^{\frac{n}{n-1}}}}{F^{\circ}(x)^{n\beta}}dx,\notag
\end{align}
which together with \eqref{FBS} implies that
\begin{equation}\label{E:3.3}
\mathfrak{T}\leq\liminf\limits_{\beta\rightarrow0}\int_{\Omega}^{} \frac{e^{\tau_{n}(1-\beta) u_{\beta}^{\frac{n}{n-1}}}}{F^{\circ}(x)^{n\beta}}dx.
\end{equation}
From the inequality $ ze^{z}\geq e^{z}-1~(z\geq0)$, we obtain
\begin{equation}
\lambda_{\beta}\geq\int_{\Omega}^{}\frac{e^{\tau_{n}(1-\beta)u_{\beta}^{\frac{n}{n-1}}}-1}{\tau_{n}(1-\beta)F^{\circ}(x)^{n\beta}}dx.\notag
\end{equation}
Then letting $\beta\to0$, one has
\begin{align}
\liminf\limits_{\beta\rightarrow 0 }\lambda_{\beta}\geq\frac{1}{\tau_{n}}\liminf\limits_{\beta\rightarrow 0 }\int_{\Omega}^{}\frac{e^{\tau_{n}(1-\beta)u_{\beta}^{\frac{n}{n-1}}}-1}{F^{\circ}(x)^{n\beta}}dx\geq\frac{1}{\tau_{n}}(\mathfrak{T}-\left |\Omega\right|).\notag
\end{align}
By \eqref{E:1.3}, we have $\mathfrak{T}>\vert\Omega\vert$. This completes the proof.  
\end{proof}

\section{Blow-up analysis} 
Let 
$$ d_{\beta}=u_{\beta}(x_{\beta})=\max_{\Omega}u_{\beta}, \quad x_{\beta}\to x_{0}\in\overline{\Omega}~\mathrm {as}~ \beta\to 0.$$
Then either 
$d_{\beta}\leq C$ or $d_{\beta}\to+\infty$ as $\beta\to0$.
If $d_{\beta}$ blows up, i.e.  $d_{\beta}\to+\infty$ as $\beta\to0$, we denote by $ \delta_{x_{0}} $ the Dirac measure centered at $x_{0}$ and 
\begin{equation}
	r_{\beta}=\lambda_{\beta}^{\frac{1}{n}}d_{\beta}^{\frac{-1}{n-1}}e^{\frac{-\tau_{n}}{n}(1-\beta)d_{\beta}^{\frac{n}{n-1}}}.\notag
\end{equation}
Then the following lemma holds:
\begin{lem}
We have $ u_{0}\equiv0 $, $ F^{n}(\nabla u_{\beta})dx\rightharpoonup\delta_{x_{0}}$ and $r_{\beta}\rightarrow0 $ as $\beta\to0$.
\end{lem} 
\begin{proof}
Let 
\begin{equation}
g_{\beta}=\frac{1}{\lambda_{\beta}}F^{\circ}(x)^{-n\beta}u_{\beta}^{\frac{1}{n-1}}e^{\tau_{n}(1-\beta)u_{\beta}^{\frac{n}{n-1}}}.\notag
\end{equation}
We know by the H$\mathrm{\ddot o}$lder inequality that
\begin{align}\label{E:3.31}
\int_{\Omega }^{}g_{\beta}^{p}dx=&\frac{1}{\lambda_{\beta}^{p}}\int_{\Omega }^{}\frac{u_{\beta}^{\frac{p}{n-1}}e^{\tau_{n}(1-\beta)pu_{\beta}^{\frac{n}{n-1}}}}{F^{\circ}(x)^{np\beta}}dx\notag\\
\leq&\frac{1}{\lambda_{\beta}^{p}}\bigg(\int_{\Omega }^{}\frac{e^{\tau_{n}(1-\beta)pp_{1}u_{\beta}^{\frac{n}{n-1}}}}{F^{\circ}(x)^{n\beta}}dx\bigg)^{1/p_{1}}\bigg(\int_{\Omega }^{}\frac{u_{\beta}^{\frac{pp_{2}}{n-1}}}{F^{\circ}(x)^{n\beta(p-\frac{1}{p_{1}})p_{2}}}dx\bigg)^{1/p_{2}},
\end{align}
where $p>1$, $p_{1}>1$ and $1/p_{1}+1/p_{2}=1$.
Assume $u_{0}\not\equiv0$. We further restrict $pp_{1}<1/(1-\Vert u_{0}\Vert_{F(\Omega)}^{n})^{1/(n-1)}$, and Lemma \ref{le:2.2} yields 
\begin{equation}\label{E:3.41}
\int_{\Omega }^{}\frac{e^{\tau_{n}(1-\beta)pp_{1}u_{\beta}^{\frac{n}{n-1}}}}{F^{\circ}(x)^{n\beta}}dx\leq C.
\end{equation}
Since $u_{\beta}\rightarrow u_{0}$ strongly in $L^{s}(\Omega)$ for $s>1$, it is obvious that
\begin{equation}\label{E:3.51}
\int_{\Omega }^{}\frac{u_{\beta}^{\frac{pp_{2}}{n-1}}}{F^{\circ}(x)^{n\beta(p-\frac{1}{p_{1}})p_{2}}}dx\leq C
\end{equation}
for $\beta$ sufficiently small. Substituting \eqref{E:3.41} and \eqref{E:3.51} into \eqref{E:3.31} and combining with Lemma \ref{le:3.1} yields $g_{\beta}\in L^{p}(\Omega)$ for some $1<p<1/(1-\Vert u_{0}\Vert_{F(\Omega)}^{n})^{1/(n-1)}$. From \eqref{E:3.1} and Lemma \ref{le:2.3}, we know that $ u_{\beta} $ is uniformly bounded
in $\Omega$. This contradicts $ d_{\beta}\rightarrow+\infty $ as $\beta\rightarrow0 $. Consequently, $ u_{0}\equiv0 $.

Now we are going to prove that $ F^{n}(\nabla u_{\beta})dx\rightharpoonup\delta_{x_{0}}$ in the sense of measure. Note that $\Vert u_{\beta}\Vert_{F(\Omega)}=1$. If $F^{n}(\nabla u_{\beta})dx\not\rightharpoonup\delta_{x_{0}}$, there exists a constant $r>0$ such that $\mathcal{W}_{r}(x_{0})\subset \Omega$ with
\begin{equation}
\limsup\limits_{\beta\rightarrow0}\Vert\zeta u_{\beta}\Vert_{F(\mathcal{W}_{r}(x_{0}))}^{n}= \iota<1,\notag 
\end{equation}
where $\zeta\in C_{0}^{1}(\Omega)$ is a cut-off function  satisfying $\zeta\equiv1$ in $\mathcal{W}_{r_{0}/2}(x_{0})$, $0\leq\zeta\leq1$ in $\mathcal{W}_{r}(x_{0})\setminus \mathcal{W}_{r_{0}/2}(x_{0})$ and $\zeta\equiv0$ in $\Omega\setminus\mathcal{W}_{r}(x_{0})$.
In view of \eqref{E:1.4}, one gets
\begin{align}
\int_{\mathcal{W}_{r}(x_{0})}^{} \frac{e^{\tau_{n}(1-\beta)\vert\frac{ \zeta u_{\beta}}{\Vert\zeta u_{\beta}\Vert_{F(\mathcal{W}_{r}(x_{0}))}}\vert^{\frac{n}{n-1}}}}{F^{\circ}(x)^{n\beta}}dx\leq&\sup\limits_{u\in W_{0}^{1,n}(\mathcal{W}_{r}(x_{0})),\Vert u\Vert_{F(\mathcal{W}_{r}(x_{0}))}\leq1}\int_{\mathcal{W}_{r}(x_{0})}^{} \frac{e^{\tau_{n}(1-\beta)\vert u\vert^{\frac{n}{n-1}}}}{F^{\circ}(x)^{n\beta}}dx\notag\\
<&+\infty,\notag
\end{align}
which implies 
\begin{equation}
\int_{\mathcal{W}_{r}(x_{0})}^{} \frac{e^{\tau_{n}q(1-\beta)\vert\zeta u_{\beta}\vert^{\frac{n}{n-1}}}}{F^{\circ}(x)^{n\beta}}dx\leq C,~\forall~1<q\leq1/\iota^{\frac{1}{n-1}}.\notag
\end{equation}
Similarly, we obtain $-\mathrm{Q}_{n} (\zeta u_{\beta})$ is bounded in $L^{q_{1}}(\mathcal{W}_{r}(x_{0}))$ for some $1<q_{1}<1/\iota^{\frac{1}{n-1}}$. By Lemma \ref{le:2.3} again, we conclude that $ u_{\beta} $ is uniformly bounded in $ \mathcal{W}_{r/2}(x_{0})$, which contradicts $ d_{\beta}\rightarrow+\infty $ as $\beta\to0$. This leads to $F^{n}(\nabla u_{\beta})dx\rightharpoonup\delta_{x_{0}}$.

Let
\begin{equation}
p_{1}>1,~0<\vartheta<\tau_{n},~\frac{1}{p_{1}}+\frac{1}{p_{2}}=1~\mathrm {and}~\vartheta p_{1}\leq\tau_{n}.\notag
\end{equation}
With the aid of \eqref{E:1.4}, one gets
\begin{align}
r_{\beta}^{n}e^{\vartheta(1-\beta)d_{\beta}^{\frac{n}{n-1}}}=&\frac{1}{d_{\beta}^{\frac{n}{n-1}}}e^{-(\tau_{n}-\vartheta)(1-\beta)d_{\beta}^{\frac{n}{n-1}}}\int_\Omega^{}\frac{u_{\beta}^{\frac{n}{n-1}}e^{\tau_{n}(1-\beta)u_{\beta}^{\frac{n}{n-1}}}}{F^{\circ}(x)^{2\beta}}dx\notag\\
\leq&\frac{1}{d_{\beta}^{\frac{n}{n-1}}}\bigg(\int_\Omega^{}\frac{e^{\vartheta p_{1}(1-\beta)u_{\beta}^{\frac{n}{n-1}}}}{F^{\circ}(x)^{n\beta}}dx\bigg)^{\frac{1}{p_{1}}}\bigg(\int_\Omega^{}\frac{u_{\beta}^{\frac{p_{2}n}{n-1}}}{F^{\circ}(x)^{n\beta}}dx\bigg)^{\frac{1}{p_{2}}}\notag\\
\leq&\frac{C}{d_{\beta}^{\frac{n}{n-1}}}\bigg(\int_\Omega^{}\frac{u_{\beta}^{\frac{p_{2}n}{n-1}}}{F^{\circ}(x)^{n\beta}}dx\bigg)^{\frac{1}{p_{2}}}.\notag 
\end{align}
Choose $q_{1}>1$ such that $1/q_{1}+1/q_{2}=1$ and $\beta q_{1}<1$. The strong convergence $ u_{\beta}\rightarrow 0 $ in $ L^{s}(\Omega) $ for every $s>1$ implies
\begin{align}
\int_\Omega^{}\frac{u_{\beta}^{\frac{p_{2}n}{n-1}}}{F^{\circ}(x)^{n\beta}}dx\leq\bigg(\int_\Omega^{}\frac{1}{F^{\circ}(x)^{n\beta q_{1}}}dx\bigg)^{\frac{1}{q_{1}}}\bigg(\int_\Omega^{}u_{\beta}^{\frac{p_{2}q_{2}n}{n-1}}dx\bigg)^{\frac{1}{q_{2}}}=o_{\beta}(1).\notag
\end{align}
Hence, we have $r_{\beta}\to0$ as $\beta \to 0$.
\end{proof}
The above lemma implies that 
\begin{align}
	&u_{\beta}\rightharpoonup 0 \quad \mathrm {weakly ~in}~ W_{0}^{1,n}(\Omega),\notag\\
	&u_{\beta}\rightarrow 0 \quad \mathrm {strongly ~in}~L^{s}(\Omega),~\forall~s>1 ,\notag\\
	&u_{\beta}\rightarrow 0 \quad\mathrm {a.e.~in}~\Omega.\notag
\end{align}

Next, we consider three cases depending on the location of the blow-up point $x_{0}$ in the region $\bar{\Omega}$. They are: 

(i)~$x_{0}\in\Omega$ and $x_{0}=0$;

(ii)~$x_{0}\in\Omega$ and $x_{0}\neq0$; 

(iii)~$x_{0}\in\partial\Omega$. \\
 For cases (i) and (ii), we are going to study the asymptotic behavior of $ u_{\beta} $ near and away from $x_{0}$ and give an upper bound estimate of $\mathfrak{T}$. And we also prove that case (iii) does not happen. 
\vspace{0.2em}
\subsection{Case (i): $x_{0}\in\Omega$ and $x_{0}=0$} In this case, either $F^{\circ}(x_{\beta})^{1-\beta}/r_{\beta}$ tends to infinity or $F^{\circ}(x_{\beta})^{1-\beta}/r_{\beta}$ is bounded.  
\subsubsection{ \textbf{If $F^{\circ}(x_{\beta})^{1-\beta}/r_{\beta}\rightarrow \infty$ as $\beta\rightarrow0$.}}
To estimate an upper bound for $\mathfrak{T}$, we divide the arguments into three steps: 

\textbf{Step 1.}~The asymptotic behavior of $ u_{\beta} $ near the point $x_{0}$.

In this step, we introduce two scaling functions for $u_{\beta}$ and prove their convergence based on the Regularity Theory, Ascoli-Arzela's Theorem and Liouville's Theorem to characterize the asymptotic behavior of $u_{\beta}$ near $x_{0}$. 
 
Let $\Omega_{\beta,1}=\{x\in \mathbb{R}^{n}:x_{\beta}+r_{\beta} F^{\circ}(x_{\beta})^{\beta}x\in \Omega\} $. We define two function sequences 
\begin{equation}
	\eta_{\beta}(x)=d_{\beta}^{-1}u_{\beta}(x_{\beta}+r_{\beta} F^{\circ}(x_{\beta})^{\beta}x)\notag 
\end{equation}
and 
\begin{equation}
\gamma_{\beta}(x)=d_{\beta}^{\frac{1}{n-1}}(u_{\beta}(x_{\beta}+r_{\beta} F^{\circ}(x_{\beta})^{\beta}x)-d_{\beta})\notag 
\end{equation}
on the set $\Omega_{\beta,1}$.
Since $F^{\circ}(x_{\beta})^{1-\beta}/r_{\beta}\to\infty$, we have $r_{\beta}F^{\circ}(x_{\beta})^{\beta}=o_{\beta}(1)F^{\circ}(x_{\beta})$, which implies that $\Omega_{\beta,1}\to\mathbb{R}^{n}$ and
\begin{eqnarray}\label{E:5.1}
\frac{F^{\circ}(x_{\beta})^{n\beta}}{F^{\circ}(x_{\beta}+r_{\beta}F^{\circ}(x_{\beta})^{\beta}x)^{n\beta}}=1+o_{\beta}(1)
\end{eqnarray}
uniformly on $\mathcal{W}_{r}$ for any fixed $R>0$ as $\beta\to0$. It is not difficult to see that $\vert\eta_{\beta}(x)\vert\leq1$ and $\gamma_{\beta}(x)\leq0$.

As can be checked by direct computation that
\begin{equation}\label{E:5.2}
-\mathrm{Q}_{n}\eta_{\beta}(x)=\frac{\eta_{\beta}^{\frac{1}{n-1}}F^{\circ}(x_{\beta})^{n\beta}}{d_{\beta}^{n}F^{\circ}(x_{\beta}+r_{\beta}F^{\circ}(x_{\beta})^{\beta}x)^{n\beta}}e^{\tau_{n}(1-\beta)(u_{\beta}^{\frac{n}{n-1}}(x_{\beta}+r_{\beta}F^{\circ}(x_{\beta})^{\beta}x)-d_{\beta}^{\frac{n}{n-1}})}~\mathrm{in}~ \Omega_{\beta,1}
\end{equation}
and
\begin{equation}
-\mathrm{Q}_{n}\gamma_{\beta}(x)=\frac{\eta_{\beta}^{\frac{1}{n-1}}F^{\circ}(x_{\beta})^{n\beta}}{F^{\circ}(x_{\beta}+r_{\beta}F^{\circ}(x_{\beta})^{\beta}x)^{n\beta}}e^{\tau_{n}(1-\beta)(u_{\beta}^{\frac{n}{n-1}}(x_{\beta}+r_{\beta}F^{\circ}(x_{\beta})^{\beta}x)-d_{\beta}^{\frac{n}{n-1}})}~ \mathrm{in}~ \Omega_{\beta,1}.\notag
\end{equation}
Note that $0\leq\eta_{\beta}(x)\leq1$, $u_{\beta}(x_{\beta}+r_{\beta}F^{\circ}(x_{\beta})^{\beta}x)\leq d_{\beta}$ and $d_{\beta}\to+\infty$ as $\beta\to0$.
It follows from \eqref{E:5.1} that $$f_{\beta,1}(x):= \frac{\eta_{\beta}^{\frac{1}{n-1}}F^{\circ}(x_{\beta})^{n\beta}}{d_{\beta}^{n}F^{\circ}(x_{\beta}+r_{\beta}F^{\circ}(x_{\beta})^{\beta}x)^{n\beta}}e^{\tau_{n}(1-\beta)(u_{\beta}^{\frac{n}{n-1}}(x_{\beta}+r_{\beta}F^{\circ}(x_{\beta})^{\beta}x)-d_{\beta}^{\frac{n}{n-1}})}\to0.$$
Therefore, we obtain $\eta_{\beta}(x)\in L^{\infty}(\mathcal{W}_{R})$ and $f_{\beta,1}\in L^{\infty}(\mathcal{W}_{R})$ for any fixed $R>0$.
Applying the Regularity Theory(see\cite{ref40}) to \eqref{E:5.2}, we conclude that $\eta_{\beta}(x)$ is uniformly bounded in $C^{1,\theta}(\mathcal{W}_{R/2})$, where $0<\theta\leq 1$. By Ascoli-Arzela's Theorem, one has 
\begin{equation}\label{E:3.77}
	\eta_{\beta}(x)\rightarrow\eta_{0}(x)~\mathrm {in}~C_{\mathrm {loc}}^{1}(\mathbb{R}^{n})~ \mathrm {as}~\beta\to0.
\end{equation}
Since  $\left|\eta_{\beta}(x)\right|\leq1$,  $\eta_{0}(0)=1$ and $ -\mathrm{Q}_{n}\eta_{0}(x)=0$ in $ \mathbb{R}^{n} $, the Liouville's Theorem(see\cite{ref41}) gives $ \eta_{0}(x)=\eta_{0}(0)\equiv1 $.

Similarly, we have 
\begin{align}
f_{\beta,2}(x):=\frac{\eta_{\beta}^{\frac{1}{n-1}}F^{\circ}(x_{\beta})^{n\beta}}{F^{\circ}(x_{\beta}+r_{\beta}F^{\circ}(x_{\beta})^{\beta}x)^{n\beta}}e^{\tau_{n}(1-\beta)(u_{\beta}^{\frac{n}{n-1}}(x_{\beta}+r_{\beta}F^{\circ}(x_{\beta})^{\beta}x)-d_{\beta}^{\frac{n}{n-1}})}\leq1+o_{\beta}(1),\notag
\end{align}
which implies $f_{\beta,2}(x)\in L^{p}(\mathcal{W}_{R})$ for some $p>1$. From Theorem 6 in \cite{ref39}, we get $\gamma_{\beta}(x)\in L^{\infty}(\mathcal{W}_{R/2})$. By Regularity Theory and Ascoli-Arzela's Theorem, we obtain
\begin{equation}\label{E:3.88}
	\gamma_{\beta}(x)\rightarrow\gamma_{0}(x)~\mathrm {in}~ C_{\mathrm {loc}}^{1}(\mathbb{R}^{n})~ \mathrm {as}~\beta\to0. 
\end{equation}  
The Mean Value Theorem leads to
\begin{align}
u_{\beta}^{\frac{n}{n-1}}(x_{\beta}+r_{\beta}F^{\circ}(x_{\beta})^{\beta}x)-d_{\beta}^{\frac{n}{n-1}}=& \frac{n}{n-1}t_{\beta}^{\frac{1}{n-1}}(u_{\beta}(x_{\beta}+r_{\beta}F^{\circ}(x_{\beta})^{\beta}x)-d_{\beta})\notag\\
=&\frac{n}{n-1}\gamma_{0}(x)+o_{\beta}(1),\notag
\end{align}
where $t_{\beta}$ lies between $u_{\beta}(x_{\beta}+r_{\beta}F^{\circ}(x_{\beta})^{\beta}x)$ and $d_{\beta}$. Therefore, we obtain 
\begin{equation}
-\mathrm{Q}_{n}\gamma_{0}(x)=e^{\frac{n}{n-1}\tau_{n}\gamma_{0}(x)}~ \mathrm {in}~ \mathbb{R}^{n},\notag
\end{equation}
where $\gamma_{0}(x)$ attains its maximum at $x=0$. 
From Lemma 4.1 of \cite{ref22}, we have 
\begin{equation}
	\gamma_{0}(x)=-\frac{n-1}{\tau_{n}}\log(1+\kappa_{n}^{\frac{1}{n-1}}F^{\circ}(x)^{\frac{n}{n-1}})\notag 
\end{equation}
and
\begin{equation}\label{E:5.4}
\int_{\mathbb{R}^{n}}^{}e^{\frac{n}{n-1}\tau_{n}\gamma_{0}(x)}dx=1.
\end{equation}
Observe that \eqref{E:3.77} and \eqref{E:3.88} imply the asymptotic behavior of $u_{\beta}$ near $x_{0}$.

\textbf{Step 2.}
The asymptotic behavior of $ u_{\beta} $ away from the point $ x_{0 }$.

In this step, we prove that $d_{\beta}^{\frac{1}{n-1}}u_{\beta}$ converges to a Green's function as $\beta\to0$, which characterizes  the asymptotic behavior of $u_{\beta}$ away from $x_{0}$.

For any $M>1$, define $$ u_{\beta,M}=\min\{\frac{d_{\beta}}{M},u_{\beta}\}. $$ 
Applying the Divergence Theorem combined and Lemma \ref{le:2.1}, we have
\begin{equation}\label{E:5.5}
-\int_{\Omega}^{}u_{\beta,M}\mathrm{Q}_{n} u_{\beta}dx=-\int_{\Omega}^{}u_{\beta,M}div(F^{n-1}(\nabla u_{\beta,M})F_{\xi}(\nabla u_{\beta,M}))=\Vert u_{\beta,M}\Vert_{F(\Omega)}^{n}
\end{equation}
and
\begin{equation}
-\int_{\Omega}^{}(u_{\beta}-\frac{d_{\beta}}{M})^{+}\mathrm{Q}_{n} u_{\beta}dx=\Vert (u_{\beta}-\frac{d_{\beta}}{M})^{+}\Vert_{F(\Omega)}^{n}.\notag
\end{equation}
Using $ u_{\beta,M} $ as a test function in \eqref{E:3.1}, we obtain
\begin{align}\label{E:5.6}
&-\int_{\Omega}^{}u_{\beta,M}\mathrm{Q}_{n} u_{\beta}dx\notag\\
&\geq\frac{1}{\lambda_{\beta}}\int_{\mathcal{W}_{Rr_{\beta}F^{\circ}(x_{\beta})^{\beta}}(x_{\beta})}^{}\frac{u_{\beta,M}u_{\beta}^{\frac{1}{n-1}}e^{\tau_{n}(1-\beta)u_{\beta}^{\frac{n}{n-1}}}}{F^{\circ}(x)^{n\beta}}dx\notag\\
&=\frac{1}{M}(1+o_{\beta}(1))\int_{\mathcal{W}_{R(0)}}^{}\frac{F^{\circ}(x_{\beta})^{n\beta}e^{\tau_{n}(1-\beta)(u_{\beta}^{\frac{n}{n-1}}(x_{\beta}+r_{\beta}F^{\circ}(x_{\beta})^{\beta}y)-d_{\beta}^{\frac{n}{n-1}})}}{F^{\circ}(x_{\beta}+r_{\beta}F^{\circ}(x_{\beta})^{\beta}y)^{n\beta}}dy\notag\\
&=\frac{1}{M}(1+o_{\beta}(1))\bigg(\int_{\mathcal{W}_{R(0)}}^{}e^{\frac{n}{n-1}\tau_{n}\gamma_{0}(y)}dy+o_{\beta,R}(1)\bigg)
\end{align}
for any $R>0$.
Testing the equation \eqref{E:3.1} with $ (u_{\beta}-\frac{d_{\beta}}{M} )^{+} $, one gets 
\begin{align}
&-\int_{\Omega}^{}(u_{\beta}-\frac{d_{\beta}}{M} )^{+}\mathrm{Q}_{n} u_{\beta}dx\notag\\
&\geq\frac{1}{\lambda_{\beta}}\int_{\mathcal{W}_{Rr_{\beta}F^{\circ}(x_{\beta})^{\beta}}(x_{\beta})}^{}\frac{(u_{\beta}-\frac{d_{\beta}}{M} )^{+}u_{\beta}^{\frac{1}{n-1}}e^{\tau_{n}(1-\beta)u_{\beta}^{\frac{n}{n-1}}}}{F^{\circ}(x)^{n\beta}}dx\notag\\
&=(1-\frac{1}{M} )(1+o_{\beta}(1))\int_{\mathcal{W}_{R(0)}}^{}\frac{F^{\circ}(x_{\beta})^{n\beta}e^{\tau_{n}(1-\beta)(u_{\beta}^{\frac{n}{n-1}}(x_{\beta}+r_{\beta}F^{\circ}(x_{\beta})^{\beta}y)-d_{\beta}^{\frac{n}{n-1}})}}{F^{\circ}(x_{\beta}+r_{\beta}F^{\circ}(x_{\beta})^{\beta}y)^{n\beta}}dy\notag\\
&=(1-\frac{1}{M} )(1+o_{\beta}(1))\bigg(\int_{\mathcal{W}_{R(0)}}^{}e^{\frac{n}{n-1}\tau_{n}\gamma_{0}(y)}dy+o_{\beta,R}(1)\bigg).\notag
\end{align}
Combining \eqref{E:5.4}, \eqref{E:5.5} and \eqref{E:5.6}, and taking the limits  $\beta\rightarrow0 $ and $ R\rightarrow+\infty $, we obtain
\begin{equation}
\liminf\limits_{\beta\rightarrow0}\Vert u_{\beta,M}\Vert_{F(\Omega)}^{n}\geq\frac{1}{M}.\notag
\end{equation}
Analogously, we have
\begin{equation}
\liminf\limits_{\beta\rightarrow0}\Vert (u_{\beta}-\frac{d_{\beta}}{M})^{+}\Vert_{F(\Omega)}^{n}\geq1-\frac{1}{M}.\notag
\end{equation}
Note that
\begin{equation}\label{E:5.7}
\Vert u_{\beta,M}\Vert_{F(\Omega)}^{n}+\Vert (u_{\beta}-\frac{d_{\beta}}{M})^{+}\Vert_{F(\Omega)}^{n}=\Vert u_{\beta}\Vert_{F(\Omega)}^{n}=1.
\end{equation}
Hence 
\begin{equation}\label{E:5.8}
\liminf\limits_{\beta\rightarrow0}\Vert u_{\beta,M}\Vert_{F(\Omega)}^{n}=\frac{1}{M}.
\end{equation}
This together with \eqref{E:1.4} implies that 
\begin{equation}\label{E:5.10}
	\int_{\Omega}^{}\frac{e^{\tau_{n}p(1-\beta)u_{\beta,M}^{\frac{n}{n-1}}}}{F^{\circ}(x)^{n\beta}}dx<+\infty,~\forall~1<p\leq M^{\frac{1}{n-1}}.
\end{equation}
 
 With the aid of \eqref{E:5.10}, we are going to prove the following lemma.
\begin{lem}\label{le:3.3}
There holds
\begin{equation}
\lim_{\beta\rightarrow0}\int_{\Omega}^{} \frac{e^{\tau_{n}(1-\beta) u_{\beta}^{\frac{n}{n-1}}}}{F^{\circ}(x)^{n\beta}}dx\leq\left|\Omega\right |+\lim_{R \to+\infty} \limsup\limits_{\beta\rightarrow0}\int_{\mathcal{W}_{Rr_{\beta}F^{\circ}(x_{\beta})^{\beta}}(x_{\beta})}^{}\frac{e^{\tau_{n}(1-\beta)u_{\beta}^{\frac{n}{n-1}}}}{F^{\circ}(x)^{n\beta}}dx.\notag
	\end{equation}
\end{lem}
\begin{proof}Given $M>1$, the domain $\Omega$ can be partitioned into:  $$\Omega=\{u_{\beta}>\frac{d_{\beta}}{M}\}\cup\{u_{\beta}\leq\frac{d_{\beta}}{M}\}.$$ 
The definition of $\lambda_{\beta}$ implies
\begin{align}\label{E:5.9}
\int_{u_{\beta}>\frac{d_{\beta}}{M}}^{}\frac{e^{\tau_{n}(1-\beta)u_{\beta}^{\frac{n}{n-1}}}}{F^{\circ}(x)^{n\beta}}dx\leq\frac{M^{\frac{n}{n-1}}}{ d_{\beta}^{\frac{n}{n-1}}}\int_{u_{\beta}>\frac{d_{\beta}}{M}}^{}\frac{u_{\beta}^{\frac{n}{n-1}}e^{\tau_{n}(1-\beta)u_{\beta}^{\frac{n}{n-1}}}}{F^{\circ}(x)^{n\beta}}dx\leq M^{\frac{n}{n-1}}\frac{\lambda_{\beta}}{d_{\beta}^{\frac{n}{n-1}}}.    
\end{align}
Note that $u_{\beta}\to0$ a.e. in $\Omega$. We have by \eqref{E:5.10}
\begin{equation}\label{E:5.11}
\int_{u_{\beta}\leq\frac{d_{\beta}}{M}}^{}\frac{e^{\tau_{n}(1-\beta)u_{\beta}^{\frac{n}{n-1}}}}{F^{\circ}(x)^{n\beta}}dx\leq\int_{\Omega}^{}\frac{e^{\tau_{n}(1-\beta)u_{\beta,M}^{\frac{n}{n-1}}}}{F^{\circ}(x)^{n\beta}}dx=\int_{\Omega}^{}\frac{1}{F^{\circ}(x)^{n\beta}}dx+o_{\beta}(1).
\end{equation}
Combining \eqref{E:5.9} and \eqref{E:5.11}, letting $ M\rightarrow1 $, $\beta\to0$, we obtain
\begin{equation}\label{E:5.12}
\lim_{\beta\rightarrow0}\int_{\Omega}^{} \frac{e^{\tau_{n}(1-\beta) u_{\beta}^{\frac{n}{n-1}}}}{F^{\circ}(x)^{n\beta}}dx\leq\left|\Omega\right |+\limsup\limits_{\beta\rightarrow0}\frac{\lambda_{\beta}}{d_{\beta}^{\frac{n}{n-1}}}. 
\end{equation}
Calculation shows
\begin{align}
\int_{\mathcal{W}_{Rr_{\beta}F^{\circ}(x_{\beta})^{\beta}}(x_{\beta})}^{}\frac{e^{\tau_{n}(1-\beta)u_{\beta}^{\frac{n}{n-1}}}}{F^{\circ}(x)^{n\beta}}dx=&\int_{\mathcal{W}_{R(0)}}^{}\frac{F^{\circ}(x_{\beta})^{n\beta}r_{\beta}^{n}e^{\tau_{n}(1-\beta)u_{\beta}^{\frac{n}{n-1}}(x_{\beta}+r_{\beta}F^{\circ}(x_{\beta})^{\beta}y)}}{F^{\circ}(x_{\beta}+r_{\beta}F^{\circ}(x_{\beta})^{\beta}y)^{n\beta}}dy\notag\\
=&\frac{\lambda_{\beta}}{d_{\beta}^{\frac{n}{n-1}}}\int_{\mathcal{W}_{R(0)}}^{}\frac{F^{\circ}(x_{\beta})^{n\beta}e^{\tau_{n}(1-\beta)(u_{\beta}^{\frac{n}{n-1}}(x_{\beta}+r_{\beta}F^{\circ}(x_{\beta})^{\beta}y)-d_{\beta}^{\frac{n}{n-1}})}}{F^{\circ}(x_{\beta}+r_{\beta}F^{\circ}(x_{\beta})^{\beta}y)^{n\beta}}dy\notag\\
=&\frac{\lambda_{\beta}}{d_{\beta}^{\frac{n}{n-1}}}(1+o_{\beta}(1))\bigg(\int_{\mathcal{W}_{R(0)}}^{}e^{\frac{n}{n-1}\tau_{n}\gamma_{0}(y)}dy+o_{\beta,R}(1)\bigg)\notag
\end{align}
for any $ R>0 $. Then using \eqref{E:5.4}, we get
\begin{align}\label{E:5.13}
\lim\limits_{R\rightarrow+\infty} \limsup\limits_{\beta\rightarrow0}\int_{\mathcal{W}_{Rr_{\beta}F^{\circ}(x_{\beta})^{\beta}}(x_{\beta})}^{}\frac{e^{\tau_{n}(1-\beta)u_{\beta}^{\frac{n}{n-1}}}}{F^{\circ}(x)^{n\beta}}dx=\limsup\limits_{\beta\rightarrow0}\frac{\lambda_{\beta}}{d_{\beta}^{\frac{n}{n-1}}}.
\end{align}
Substituting \eqref{E:5.13} into \eqref{E:5.12} yields Lemma \ref{le:3.3}. Moreover,
we have
\begin{equation}\label{le:3.4}
	\lim\limits_{\beta\rightarrow0}\frac{\lambda_{\beta}}{d_{\beta}}=+\infty~\mathrm {and}~\limsup\limits_{\beta\rightarrow0}\frac{d_{\beta}^{\frac{n}{n-1}}}{\lambda_{\beta}}\leq C.
	\end{equation}
In fact, we argue by contradiction to suppose that $ \lambda_{\beta}/d_{\beta}^{\frac{n}{n-1}}\to0$ as $\beta\to0$. It follows from \eqref{E:3.3} and \eqref{E:5.12} that
\begin{align}
\mathfrak{T}\leq\liminf\limits_{\beta\rightarrow0}\int_{\Omega}^{} \frac{e^{\tau_{n}(1-\beta) u_{\beta}^{\frac{n}{n-1}}}}{F^{\circ}(x)^{n\beta}}dx\leq\left|\Omega\right|.\notag
\end{align}
This contradicts the fact that  $\vert\Omega\vert<\mathfrak{T}<+\infty$. Consequently, we obtain $\limsup\limits_{\beta\rightarrow0}\frac{d_{\beta}^{\frac{n}{n-1}}}{\lambda_{\beta}}\leq C$.
As a byproduct, it is easy to know $\lim\limits_{\beta\rightarrow0}\frac{\lambda_{\beta}}{d_{\beta}}=+\infty$.
\end{proof}

Next, we investigate the convergence of $ d_{\beta}^{\frac{1}{n-1}}u_{\beta}$ as $\beta\to0$.
\begin{lem}\label{le:3.5}
The family $ d_{\beta}^{\frac{1}{n-1}}u_{\beta}\rightharpoonup G_{0}$ weakly in $W_{0}^{1,q}(\Omega)$ for any $ 1<q<n $ and $ d_{\beta}^{\frac{1}{n-1}}u_{\beta}\rightarrow G_{0}$ in $ C_{\mathrm{loc}}^{1}(\overline{\Omega}\setminus\{0\}) $. Here, $G_{0}$ is a Green's function solving $$ -\mathrm{Q}_{n}G_{0}=\delta_{0}$$ in a distributional sense. Moreover, $G_{0}$ can be written as 
\begin{equation}\label{E:5.14}
G_{0}(x)=-\frac{n}{\tau_{n}}\log F^{\circ}(x)+C_{0}+D(x),  
\end{equation} 
where $C_{0}$ is a constant, $ D(x)\in C^{0}(\overline{\Omega})\cap C_{\mathrm{loc}}^{1}(\overline{\Omega}\setminus\{0\})$ satisfies $D(0)=0 $ and $$\lim\limits_{x\rightarrow0}F^{\circ}(x)\nabla D(x)=0.$$
\end{lem}
\begin{proof}
Given any $\varpi(x)\in C^{\infty}(\overline{\Omega})$, we shall show that
\begin{equation}\label{E:5.15}
\lim\limits_{\beta\rightarrow0}\int_{\Omega}^{}\frac{d_{\beta}u_{\beta}^{\frac{1}{n-1}}e^{\tau_{n}(1-\beta)u_{\beta}^{\frac{n}{n-1}}}}{\lambda_{\beta}F^{\circ}(x)^{n\beta}}\varpi(x)dx=\varpi(0).
\end{equation}
The limiting behavior $F^{\circ}(x_{\beta})^{1-\beta}/r_{\beta}\to\infty$ as $\beta\to0$ implies that $ \mathcal{W}_{Rr_{\beta}F^{\circ}(x_{\beta})^{\beta}}(x_{\beta})\subset\{u_{\beta}>\frac{d_{\beta}}{M}\}$ for $ \beta$ sufficiently close to 0. We denote
\begin{align}
\textup{I}=&\int_{u_{\beta}\leq\frac{d_{\beta}}{M}}^{}\frac{d_{\beta}u_{\beta}^{\frac{1}{n-1}}e^{\tau_{n}(1-\beta)u_{\beta}^{\frac{n}{n-1}}}}{\lambda_{\beta}F^{\circ}(x)^{n\beta}}\varpi(x)dx,\notag\\
\textup{II}=&\int_{\mathcal{W}_{Rr_{\beta}F^{\circ}(x_{\beta})^{\beta}}(x_{\beta})}^{}\frac{d_{\beta}u_{\beta}^{\frac{1}{n-1}}e^{\tau_{n}(1-\beta)u_{\beta}^{\frac{n}{n-1}}}}{\lambda_{\beta}F^{\circ}(x)^{n\beta}}\varpi(x)dx\notag
\end{align}
and 
\begin{equation}
\textup{III}=\int_{\{u_{\beta}>\frac{d_{\beta}}{M}\}\setminus  \mathcal{W}_{Rr_{\beta}F^{\circ}(x_{\beta})^{\beta}}(x_{\beta})}^{}\frac{d_{\beta}u_{\beta}^{\frac{1}{n-1}}e^{\tau_{n}(1-\beta)u_{\beta}^{\frac{n}{n-1}}}}{\lambda_{\beta}F^{\circ}(x)^{n\beta}}\varpi(x)dx.\notag
\end{equation}
A simple estimate yields
\begin{align}
\textup{I}\leq&\frac{d_{\beta}}{\lambda_{\beta}}\left(\sup\limits_{\overline{\Omega}}\left|\varpi(x)\right|\right)\int_{u_{\beta}\leq\frac{d_{\beta}}{M}}^{}\frac{u_{\beta}^{\frac{1}{n-1}}e^{\tau_{n}(1-\beta)u_{\beta}^{\frac{n}{n-1}}}}{F^{\circ}(x)^{n\beta}}dx\notag\\
\leq&\frac{d_{\beta}}{\lambda_{\beta}}\left(\sup\limits_{\overline{\Omega}}\left|\varpi(x)\right|\right)\bigg(\int_{\Omega}^{}\frac{e^{\tau_{n}(1-\beta)p_{1}u_{\beta,M}^{\frac{n}{n-1}}}}{F^{\circ}(x)^{n\beta}}dx\bigg)^{1/p_{1}}\bigg(\int_{\Omega}^{}\frac{u_{\beta,M}^{\frac{p_{2}}{n-1}}}{F^{\circ}(x)^{n\beta}}dx\bigg)^{1/p_{2}}.\notag
\end{align}
Choosing $p_{1}\leq M^{\frac{1}{n-1}}$ such that $1/p_{1}+1/p_{2}=1$, we deduce from \eqref{E:5.10} and
 \eqref{le:3.4} that
\begin{equation}\label{E:5.16}
\lim\limits_{\beta\rightarrow0}\textup{I}=0.
\end{equation}
A straightforward calculation shows that
\begin{align}
\textup{II}=&(1+o_{\beta}(1))\int_{\mathcal{W}_{R(0)}}^{}\frac{\varpi(x_{\beta}+r_{\beta}F^{\circ}(x_{\beta})^{\beta}y)F^{\circ}(x_{\beta})^{n\beta}r_{\beta}^{n}d_{\beta}^{\frac{n}{n-1}}e^{\tau_{n}(1-\beta)u_{\beta}^{\frac{n}{n-1}}(x_{\beta}+r_{\beta}F^{\circ}(x_{\beta})^{\beta}y)}}{\lambda_{\beta}F^{\circ}(x_{\beta}+r_{\beta}F^{\circ}(x_{\beta})^{\beta}y)^{n\beta}}dy\notag\\
=&(1+o_{\beta}(1))\int_{\mathcal{W}_{R(0)}}^{}\frac{\varpi(x_{\beta}+r_{\beta}F^{\circ}(x_{\beta})^{\beta}y)F^{\circ}(x_{\beta})^{n\beta}e^{\tau_{n}(1-\beta)(u_{\beta}^{\frac{n}{n-1}}(x_{\beta}+r_{\beta}F^{\circ}(x_{\beta})^{\beta}y)-d_{\beta}^{\frac{n}{n-1}})}}{F^{\circ}(x_{\beta}+r_{\beta}F^{\circ}(x_{\beta})^{\beta}y)^{n\beta}}dy\notag\\
=&\varpi(x_{\beta}+r_{\beta}F^{\circ}(x_{\beta})^{\beta}y)(1+o_{\beta}(1))\int_{\mathcal{W}_{R(0)}}^{}e^{\tau_{n}(1-\beta)(u_{\beta}^{\frac{n}{n-1}}(x_{\beta}+r_{\beta}F^{\circ}(x_{\beta})^{\beta}y)-d_{\beta}^{\frac{n}{n-1}})}dy\notag\\
=&\varpi(x_{\beta}+r_{\beta}F^{\circ}(x_{\beta})^{\beta}y)(1+o_{\beta}(1))\bigg(\int_{\mathcal{W}_{R(0)}}^{}e^{\frac{n}{n-1}\tau_{n}\gamma_{0}(y)}dy+o_{\beta,R}(1)\bigg).\notag
\end{align}
Passing to the limit as $ \beta\rightarrow0 $ and $ R\rightarrow+\infty $,  we obtain
\begin{equation}\label{E:5.17}
\lim\limits_{\beta\rightarrow0}\textup{II}=\varpi(0).
\end{equation}
In addition, one sees that
\begin{align}
\textup{III}\leq&M\int_{\{u_{\beta}>\frac{d_{\beta}}{M}\}\setminus \mathcal{W}_{Rr_{\beta}F^{\circ}(x_{\beta})^{\beta}}(x_{\beta})}^{}\frac{u_{\beta}^{\frac{n}{n-1}}e^{\tau_{n}(1-\beta)u_{\beta}^{\frac{n}{n-1}}}}{\lambda_{\beta}F^{\circ}(x)^{n\beta}}\varpi(x)dx\notag\\
\leq&M\bigg(\sup\limits_{\overline{\Omega}}\left|\varpi(x)\right|\bigg)\lambda_{\beta}^{-1}\bigg(\int_{\Omega}^{}\frac{u_{\beta}^{\frac{n}{n-1}}e^{\tau_{n}(1-\beta)u_{\beta}^{\frac{n}{n-1}}}}{F^{\circ}(x)^{n\beta}}dx-\int_{\mathcal{W}_{Rr_{\beta}F^{\circ}(x_{\beta})^{\beta}}(x_{\beta})}^{}\frac{u_{\beta}^{\frac{n}{n-1}}e^{\tau_{n}(1-\beta)u_{\beta}^{\frac{n}{n-1}}}}{F^{\circ}(x)^{n\beta}}dx\bigg)\notag\\
=&M\bigg(\sup\limits_{\overline{\Omega}}\left|\varpi(x)\right|\bigg)\bigg(1-\int_{\mathcal{W}_{R(0)} }^{}e^{\frac{n}{n-1}\tau_{n}\gamma_{0}(y)}dy+o_{\beta,R}(1)\bigg).\notag
\end{align}
Similarly, we have $\lim\limits_{\beta\rightarrow0}\textup{III}=0$, which together with \eqref{E:5.16} and \eqref{E:5.17} gives \eqref{E:5.15}. 

By \eqref{E:3.1}, one gets
\begin{equation}
-\mathrm{Q}_{n}(d_{\beta}^{\frac{1}{n-1}}u_{\beta})=\frac{d_{\beta}u_{\beta}^{\frac{1}{n-1}}e^{\tau_{n}(1-\beta)u_{\beta}^{\frac{n}{n-1}}}}{\lambda_{\beta}F^{\circ}(x)^{n\beta}}.\notag
\end{equation}
Choosing $\varpi(x)\equiv1$ in \eqref{E:5.15}, we find that $\mathrm{Q}_{n}(d_{\beta}^{\frac{1}{n-1}}u_{\beta})$ is bounded in $ L^{1}(\Omega)$. From Lemma \ref{le:2.4}, one can see that $d_{\beta}^{\frac{1}{n-1}}u_{\beta}$ is bounded in $W_{0}^{1,q}(\Omega)$ for $ 1<q<n $. Then
we have $d_{\beta}^{\frac{1}{n-1}}u_{\beta}\rightharpoonup G_{0}$ weakly in $W_{0}^{1,q}(\Omega)$ and $d_{\beta}^{\frac{1}{n-1}}u_{\beta}\rightarrow G_{0}$ in $ C_{\mathrm{loc}}^{1}(\overline{\Omega}\setminus\{0\}) $. Here $ G_{0} $ satisfies 
\begin{equation}
	\left\{  
	\begin{aligned}
&-\mathrm{Q}_{n} G_{0}=\delta_{0}~\mathrm {in}~\Omega,\\
&G_{0}=0~\mathrm {on}~\partial\Omega.
	\end{aligned}
	\right.\notag
\end{equation}
From \cite{ref20}, $G_{0}$ can be written as  \eqref{E:5.14}. This finishes the proof.
\end{proof}
Up to now, we had gotten the asymptotic behavior of $u_{\beta}$ away from or near $x_{0}$. 

\textbf{Step 3.} We establish an upper bound for $\mathfrak{T}$.
\begin{lem}\label{le:3.6}
Let $\mu_{\beta}\in W_{0}^{1,n}(\mathcal{W}_{1})$ satisfy
$\Vert\mu_{ \beta}\Vert_{F(\mathcal{W}_{1})}\leq1$ and $\mu_{ \beta}\rightharpoonup0$ weakly in $W_{0}^{1,n}(\mathcal{W}_{1})$ as $\beta\rightarrow0$. There holds
\begin{equation}
\limsup\limits_{\beta\rightarrow0}\int_{\mathcal{W}_{1}}^{}\frac{e^{\tau_{n}(1-\beta)\vert\mu_{\beta}\vert^{\frac{n}{n-1}}}-1}{F^{\circ}(x)^{n\beta}} dx\leq \kappa_{n}e^{{\sum_{j=1}^{n-1}}\frac{1}{j}}.\notag
	\end{equation}
\end{lem}
\begin{proof}
Given the definitions \eqref{E:2.2} and \eqref{E:2.1} for $\mu_{\beta}^{\sharp}$ and $\mu_{\beta}^{\star}$, we have $\mu_{ \beta}^{\star}\in W_{0}^{1,n}(\mathcal{W}_{1})$ and $$\Vert\mu_{ \beta}^{\star}\Vert_{F(\mathcal{W}_{1})}\leq\Vert\mu_{ \beta}\Vert_{F(\mathcal{W}_{1})}.$$
Denote $a=(\kappa_{n}/V_{n})^{1/n}$. As shown in \cite{ref23}, we obtain $\mu_{\beta}^{\sharp}\in W_{0}^{1,n}(\mathbb{B}_{a})$,
\begin{equation}\label{E:5.18}
\Vert \mu_{\beta}^{\star}\Vert_{F(\mathcal{W}_{1})}=a\Vert \mu_{\beta}^{\sharp}\Vert_{W_{0}^{1,n}(\mathbb{B}_{a})},
\end{equation}
\begin{equation}\label{E:5.19}
\int_{\mathcal{W}_{1}}^{}\frac{e^{\tau_{n}(1-\beta)\vert\mu_{\beta}^{\star}\vert^{\frac{n}{n-1}}}-1}{F^{\circ}(x)^{n\beta}}dx=a^{n\beta}\int_{\mathbb{B}_{a}}^{}\frac{e^{\tau_{n}(1-\beta)\vert\mu_{\beta}^{\sharp}\vert^{\frac{n}{n-1}}}-1}{\vert x\vert^{n\beta}}dx
\end{equation}
and 
\begin{equation}\label{E:5.20}
\limsup\limits_{\beta\rightarrow0}\int_{\mathcal{W}_{1}}^{}\frac{e^{\tau_{n}(1-\beta)\vert\mu_{\beta}\vert^{\frac{n}{n-1}}}-1}{F^{\circ}(x)^{n\beta}}dx=\limsup\limits_{\beta\rightarrow0}\int_{\mathcal{W}_{1}}^{}\frac{e^{\tau_{n}(1-\beta)\vert\mu_{\beta}^{\star}\vert^{\frac{n}{n-1}}}-1}{F^{\circ}(x)^{n\beta}}dx.
	\end{equation}
The scaling $v_{k}(x)=a\mu_{k}(ax)$ leads to
\begin{equation}
\Vert v_{\beta}^{\sharp}\Vert_{W_{0}^{1,n}(\mathbb{B})}=a\Vert \mu_{\beta}^{\sharp}\Vert_{W_{0}^{1,n}(\mathbb{B}_{a})}.\notag
\end{equation}
This together with \eqref{E:5.18} leads to
\begin{equation}
\Vert \mu_{\beta}^{\star}\Vert_{F(\mathcal{W}_{1})}=\Vert v_{\beta}^{\sharp}\Vert_{W_{0}^{1,n}(\mathbb{B})}.\notag
\end{equation}
A direct calculation shows
\begin{align}\label{E:5.21}
\int_{\mathbb{B}}^{}\frac{e^{\alpha_{n}(1-\beta)\vert v_{\beta}^{\sharp}\vert^{\frac{n}{n-1}}}-1}{\vert x\vert^{n\beta}}dx&=nV_{n}\int_{0}^{1}\frac{e^{\alpha_{n}(1-\beta)\vert a\mu_{\beta}^{\sharp}(ar)\vert^{\frac{n}{n-1}}}-1}{\vert r\vert^{n\beta}} r^{n-1}dr\notag\\
&=a^{n(\beta-1)}nV_{n}\int_{0}^{a}\frac{e^{\tau_{n}(1-\beta)\vert\mu_{\beta}^{\sharp}(y)\vert^{\frac{n}{n-1}}}-1}{\vert y\vert^{n\beta}}y^{n-1}dy\notag\\
&=a^{n(\beta-1)}\int_{\mathbb{B}_{a}}^{}\frac{e^{\tau_{n}(1-\beta)\vert\mu_{\beta}^{\sharp}(x)\vert^{\frac{n}{n-1}}}-1}{\vert x\vert^{n\beta}}dx.
\end{align}
Combining \eqref{E:5.19} and \eqref{E:5.21}, we arrive at
\begin{equation}
\int_{\mathcal{W}_{1}}^{}\frac{e^{\tau_{n}(1-\beta)\vert\mu_{\beta}^{\star}\vert^{\frac{n}{n-1}}}-1}{F^{\circ}(x)^{n\beta}}dx=\frac{\kappa_{n}}{V_{n}}\int_{\mathbb{B}}^{}\frac{e^{\alpha_{n}(1-\beta)\vert v_{\beta}^{\sharp}\vert^{\frac{n}{n-1}}}-1}{\vert x\vert^{n\beta}}dx.\notag
\end{equation}
It follows from \eqref{E:5.20} and Lemma 10 in \cite{ref16} that
\begin{align}
\limsup\limits_{\beta\rightarrow0}\int_{\mathcal{W}_{1}}^{}\frac{e^{\tau_{n}(1-\beta)\vert\mu_{\beta}\vert^{\frac{n}{n-1}}}-1}{F^{\circ}(x)^{n\beta}}dx=\frac{\kappa_{n}}{V_{n}}\limsup\limits_{\beta\rightarrow0}\int_{\mathbb{B}}^{}\frac{e^{\alpha_{n}(1-\beta)\vert v_{k}^{\sharp}\vert^{\frac{n}{n-1}}}-1}{\vert x\vert^{n\beta}}dx\leq\kappa_{n}e^{{\sum_{j=1}^{n-1}}\frac{1}{j}},\notag
\end{align}
as desired.
\end{proof}

Now we are going to estimate the upper bounded of $\mathfrak{T}$.
\begin{lem}\label{le:3.7}
Let $C_{0}$ be given by \eqref{E:5.14}, then we have
\begin{equation}
\mathfrak{T}\leq\vert\Omega\vert+\kappa_{n}e^{\tau_{n}C_{0}+{\sum_{j=1}^{n-1}}\frac{1}{j}}.\notag
	\end{equation}
\end{lem}
\begin{proof}Define $u_{\beta}(\rho)=\sup_{\partial \mathcal{W}_{\rho}}u_{\beta} $ and $\tilde{u}_{\beta}=(u_{\beta}-u_{\beta}(\rho))^{+}$. Obviously,
\begin{equation}
u_{\beta}(\rho)=\frac{1}{d_{\beta}^{\frac{1}{n-1}}}(\sup_{\partial \mathcal{W}_{\rho}}G_{0}(\rho)+o_{\beta}(1)),~\tilde{u}_{\beta}\in
W_{0}^{1,n}(\mathcal{W}_{\rho})~ \mathrm {and}~\Vert\tilde{u}_{\beta}\Vert_{F(\mathcal{W}_{\rho})}=\Vert u_{\beta}\Vert_{_{F(\mathcal{W}_{\rho})}}.\notag
\end{equation} 
Denoting by $\vec{n}$ the unit outward normal to $\partial \mathcal{W}_{\rho}$, we apply Lemma \ref{le:3.5} and the divergence theorem to obtain
\begin{align}\label{E:5.22}
\Vert u_{\beta}\Vert_{F(\Omega\setminus \mathcal{W}_{\rho})}^{n}=&\frac{1}{d_{\beta}^{\frac{n}{n-1}}}\int_{\Omega\setminus \mathcal{W}_{\rho}}^{}\left(F^{n}(\nabla G_{0})+o_{\beta}(1)\right)dx\notag\\
=&\frac{1}{d_{\beta}^{\frac{n}{n-1}}}\bigg(-\int_{\Omega\setminus \mathcal{W}_{\rho}}^{}G_{0}\mathrm{Q}_{n}G_{0}dx-\int_{\partial \mathcal{W}_{\rho}}^{}G_{0}F^{^{n-1}}(\nabla G_{0})\left \langle F_{\xi}(\nabla G_{0}),\vec{n} \right\rangle ds\bigg)\notag\\
&+\frac{o_{\beta,\rho}(1)}{d_{\beta}^{\frac{n}{n-1}}}.
\end{align}
Since $ G_{0} $ is a distributional solution of $ -\mathrm{Q}_{n} G_{0}=\delta_{0}$, we obtain
\begin{equation}\label{E:5.23}
	-\int_{\Omega\setminus \mathcal{W}_{\rho}}^{}G_{0}\mathrm{Q}_{n}G_{0}dx=o_{\rho}(1).
\end{equation}
By \eqref{E:5.14} and Lemma \ref{le:2.1}, we have on $\partial \mathcal{W}_{\rho}$ 
\begin{align}
F^{^{n-1}}(\nabla G_{0})=\frac{1}{n\kappa_{n}\rho^{n-1}}+o_{\rho}(1)\notag
\end{align}
and
\begin{align}
\left\langle F_{\xi}(\nabla G_{0}),\vec{n} \right\rangle=-\frac{1}{\vert\nabla F^{\circ}(x) \vert}+o_{\rho}(1).\notag
\end{align}
The coarea
formula with respect to $F^{\circ}(x)$ yields
\begin{align}\label{E:5.24}
-\int_{\partial \mathcal{W}_{\rho}}^{}G_{0}F^{n-1}(\nabla G_{0})\left \langle F_{\xi}(\nabla G_{0}),\vec{n} \right\rangle ds&=\frac{G_{0}(\rho)}{n\kappa_{n}\rho^{n-1}}\int_{\partial \mathcal{W}_{\rho}}^{}\frac{1}{\vert\nabla F^{\circ}(x) \vert}ds+o_{\rho}(1)\notag\\
&=-\frac{n}{\tau_{n}}\log \rho+C_{0}+o_{\rho}(1).
\end{align}
In view of \eqref{E:5.22}, \eqref{E:5.23} and \eqref{E:5.24}, we obtain
\begin{equation}
\Vert u_{\beta}\Vert_{F(\Omega\setminus \mathcal{W}_{\rho})}^{n}=\frac{1}{d_{\beta}^{\frac{n}{n-1}}}\bigg(-\frac{n}{\tau_{n}}\log \rho+C_{0}+o_{\beta,\rho}(1)\bigg),\notag
\end{equation}
which implies
\begin{align}
\Vert \tilde{u}_{\beta}\Vert_{F(\mathcal{W}_{\rho})}=&(\Vert u_{\beta}\Vert_{F(\Omega)}^{n}-\Vert u_{\beta}\Vert_{F(\Omega\setminus \mathcal{W}_{\rho})}^{n})^{1/n}\notag\\
=&\bigg(1-\frac{1}{d_{\beta}^{\frac{n}{n-1}}}\bigg(-\frac{n}{\tau_{n}}\log \rho+C_{0}+o_{\beta,\rho}(1)\bigg)\bigg)^{1/n}.\notag
\end{align}
We write $s_{\beta}=\Vert \tilde{u}_{\beta}\Vert_{F(\mathcal{W}_{\rho})}$. Since  $r_{\beta}F^{\circ}(x_{\beta})^{\beta}=o_{\beta}(1)F^{\circ}(x_{\beta})$, we get $\mathcal{W}_{Rr_{\beta}F^{\circ}(x_{\beta})^{\beta}}(x_{\beta})\subset \mathcal{W}_{\rho}$ for $\beta$ sufficiently small. In $\mathcal{W}_{Rr_{\beta}F^{\circ}(x_{\beta})^{\beta}}(x_{\beta})$, we obtain
\begin{align}
\tau_{n}(1-\beta) u_{\beta}^{\frac{n}{n-1}}\leq&\tau_{n}(1-\beta) (\tilde{u}_{\beta}+u_{\beta}(\rho))^{\frac{n}{n-1}}\notag\\
\leq&\tau_{n}(1-\beta)\bigg(\tilde{u}_{\beta}^{\frac{n}{n-1}}+\frac{n}{n-1}\tilde{u}_{\beta}^{\frac{1}{n-1}}u_{\beta}(\rho)+o_{\beta}(1)\bigg) \notag\\
\leq&\tau_{n}(1-\beta)\bigg(\frac{\tilde{u}_{\beta}}{s_{\beta}}\bigg)^{\frac{n}{n-1}}-\frac{\tau_{n}(1-\beta)}{n-1}\bigg(-\frac{n}{\tau_{n}}\log \rho+C_{0}+o_{\beta,\rho}(1)\bigg)\notag\\
&+\frac{n\tau_{n}(1-\beta)}{n-1}\bigg(-\frac{n}{\tau_{n}}\log \rho+C_{0}+o_{\beta,\rho}(1)\bigg)\notag\\
=&\tau_{n}(1-\beta)\bigg(\bigg(\frac{\tilde{u}_{\beta}}{s_{\beta}}\bigg)^{\frac{n}{n-1}}-\frac{n}{\tau_{n}}\log\rho+C_{0}+o_{\beta,\rho}(1)\bigg).\notag
\end{align}
Here we used the inequality  $(1-z)^{t}\leq1-tz$ for $t,z\in(0,1)$. Therefore,
\begin{align}\label{E:5.26}
&\int_{\mathcal{W}_{Rr_{\beta}F^{\circ}(x_{\beta})^{\beta}}(x_{\beta})}^{}\frac{e^{\tau_{n}(1-\beta)u_{\beta}^{\frac{n}{n-1}}}}{F^{\circ}(x)^{n\beta}}dx\notag\\
\leq&\int_{\mathcal{W}_{Rr_{\beta}F^{\circ}(x_{\beta})^{\beta}}(x_{\beta})}^{}\frac{e^{\tau_{n}(1-\beta)((\tilde{u}_{\beta}/s_{\beta})^{\frac{n}{n-1}}-\frac{n}{\tau_{n}}\log\rho+C_{0}+o_{\beta,\rho}(1))}}{F^{\circ}(x)^{n\beta}}dx\notag\\
\leq&\rho^{-n(1-\beta)}e^{\tau_{n}(1-\beta) C_{0}+o_{\beta,\rho}(1)}\int_{\mathcal{W}_{Rr_{\beta}F^{\circ}(x_{\beta})^{\beta}}(x_{\beta})}^{}\frac{e^{\tau_{n}(1-\beta)(\tilde{u}_{\beta}/s_{\beta})^{\frac{n}{n-1}}}}{F^{\circ}(x)^{n\beta}}dx
\end{align}
It follows by Lemma \ref{le:3.6} that
\begin{align}\label{E:5.27}
\int_{\mathcal{W}_{Rr_{\beta}F^{\circ}(x_{\beta})^{\beta}}(x_{\beta})}^{}\frac{e^{\tau_{n}(1-\beta)(\tilde{u}_{\beta}/s_{\beta})^{\frac{n}{n-1}}}}{F^{\circ}(x)^{n\beta}}dx
&\leq\int_{\mathcal{W}_{\rho}}^{}\frac{e^{\tau_{n}(1-\beta)(\tilde{u}_{\beta}/s_{\beta})^{\frac{n}{n-1}}}-1}{F^{\circ}(x)^{n\beta}}dx+o_{\beta}(1)\notag\\
&\leq \rho^{n}\kappa_{n}e^{{\sum_{j=1}^{n-1}}\frac{1}{j}}.
\end{align}
By substituting \eqref{E:5.27} into \eqref{E:5.26} with $\beta\to0 $ and  $R\to+\infty$, we derive
\begin{equation}\label{E:5.28}
\lim_{R\to+\infty}\limsup\limits_{\beta\rightarrow0}\int_{\mathcal{W}_{Rr_{\beta}F^{\circ}(x_{\beta})^{\beta}}(x_{\beta})}^{}\frac{e^{\tau_{n}(1-\beta)u_{\beta}^{\frac{n}{n-1}}}}{F^{\circ}(x)^{n\beta}}dx\leq\kappa_{n}e^{\tau_{n}C_{0}+{\sum_{j=1}^{n-1}}\frac{1}{j}}. 
\end{equation}
Combining  \eqref{E:3.3}, \eqref{E:5.28} and Lemma \ref{le:3.3}, we prove Lemma \ref{le:3.7}.
\end{proof}
\vspace{0.2em}
\subsubsection{\textbf{$F^{\circ}(x_{\beta})^{1-\beta}/r_{\beta}$ is bounded as $\beta\rightarrow0$.}} We repeat the arguments in subsubsection 3.1.1.

\textbf{Step 1.} The asymptotic behavior of $u_{\beta}$ near $x_{0}$.

We define the functions
\begin{equation}
\phi_{\beta}(x)=d_{\beta}^{-1}u_{\beta}(x_{\beta}+r_{\beta}^{\frac{1}{1-\beta}}x)~\mathrm {and}~\varphi_{\beta}(x)=d_{\beta}^{\frac{1}{n-1}}(u_{\beta}(x_{\beta}+r_{\beta}^{\frac{1}{1-\beta}}x)-d_{\beta})\notag
\end{equation}
in $ \Omega_{\beta,2}=\{x\in\mathbb{R}^{n}:x_{\beta}+r_{\beta}^{\frac{1}{1-\beta}}x\in\Omega\}$. Clearly,
\begin{equation}
	\vert\phi_{\beta}(x)\vert\leq1,~ \varphi_{\beta}(x)\leq0~\mathrm {and}~\Omega_{\beta,2}\to\mathbb{R}^{n}~\mathrm {as}~ \beta\to 0.\notag
\end{equation}
It can be easily checked that
\begin{equation}\label{E:5.29}
-\mathrm{Q}_{n}\phi_{\beta}(x)=d_{\beta}^{-n}F^{\circ}\bigg(x+\frac{x_{\beta}}{r_{\beta}^{\frac{1}{1-\beta}}}\bigg)^{-n\beta}\phi_{\beta}^{\frac{1}{n-1}}e^{\tau_{n}(1-\beta)(u_{\beta}^{\frac{n}{n-1}}(x_{\beta}+r_{\beta}^{\frac{1}{1-\beta}}x)-d_{\beta}^{\frac{n}{n-1}})}~ \mathrm{in}~ \Omega_{\beta,2},
\end{equation}
and
\begin{equation}\label{E:5.30}
-\mathrm{Q}_{n}\varphi_{\beta}(x)=F^{\circ}\bigg(x+\frac{x_{\beta}}{r_{\beta}^{\frac{1}{1-\beta}}}\bigg)^{-n\beta}\phi_{\beta}^{\frac{1}{n-1}}e^{\tau_{n}(1-\beta)(u_{\beta}^{\frac{n}{n-1}}(x_{\beta}+r_{\beta}^{\frac{1}{1-\beta}}x)-d_{\beta}^{\frac{n}{n-1}})}~ \mathrm{in}~ \Omega_{\beta,2}. 
\end{equation}
Since $F^{\circ}(x_{\beta})^{1-\beta}/r_{\beta}\leq C$ as $\beta\to0$, we have $ F^{\circ}\bigg(\frac{x_{\beta}}{r_{\beta}^{\frac{1}{1-\beta}}}\bigg)\leq C$, which implies that $F^{\circ}\bigg(x+\frac{x_{\beta}}{r_{\beta}^{\frac{1}{1-\beta}}}\bigg)^{-n\beta}$ is bounded in $L^{q}(\mathcal{W}_{R})$ for any $q>1$ and fixed $R>0$.  Recalling $\left|\phi_{\beta}(x)\right|\leq1$, $u_{\beta}^{\frac{n}{n-1}}(x_{\beta}+r_{\beta}^{\frac{1}{1-\beta}}x)\leq d_{\beta}^{\frac{n}{n-1}} $ and $ d_{\beta}\rightarrow+\infty $ as $\beta\to0$, we have
\begin{equation}
k_{\beta,1}(x):=d_{\beta}^{-n}F^{\circ}\bigg(x+\frac{x_{\beta}}{r_{\beta}^{\frac{1}{1-\beta}}}\bigg)^{-n\beta}\phi_{\beta}^{\frac{1}{n-1}}e^{\tau_{n}(1-\beta)(u_{\beta}^{\frac{n}{n-1}}(x_{\beta}+r_{\beta}^{\frac{1}{1-\beta}}x)-d_{\beta}^{\frac{n}{n-1}})}\to0.\notag
\end{equation} 
Applying the Regularity Theory to \eqref{E:5.29}, one obtains 
\begin{equation}
	\phi_{\beta}(x)\rightarrow\phi_{0}(x)~\mathrm {in}~C_{\mathrm {loc}}^{1}(\mathbb{R}^{n})~\mathrm {as}~\beta\to0,\notag
\end{equation} 
where $\phi_{0}(x)$ satisfies 
$-\mathrm{Q}_{n}\phi_{0}(x)=0$.
The Liouville's Theorem implies that $\phi_{0}(x)=\phi_{0}(0)\equiv1 $.

On the other hand, our estimate shows that
\begin{equation}
k_{\beta,2}(x):=F^{\circ}\bigg(x+\frac{x_{\beta}}{r_{\beta}^{\frac{1}{1-\beta}}}\bigg)^{-n\beta}\phi_{\beta}^{\frac{1}{n-1}}e^{\tau_{n}(1-\beta)(u_{\beta}^{\frac{n}{n-1}}(x_{\beta}+r_{\beta}^{\frac{1}{1-\beta}}x)-d_{\beta}^{\frac{n}{n-1}})}\leq F^{\circ}\bigg(x+\frac{x_{\beta}}{r_{\beta}^{\frac{1}{1-\beta}}}\bigg)^{-n\beta}.\notag
\end{equation}
Clearly, we have $k_{\beta,2}(x)\in L^{q}(\mathcal{W}_{R})$ for some $q>1$ and any fixed $R>0$. Applying the Regularity Theory to \eqref{E:5.30}, we get
\begin{equation}
\varphi_{\beta}(x)\rightarrow\varphi_{0}(x)~\mathrm {in}~C_{\mathrm {loc}}^{1}(\mathbb{R}^{n})~\mathrm {as}~\beta\to0,\notag
\end{equation} 
where $\varphi_{0}(x) $ satisfies
\begin{equation}
	\left\{  
\begin{aligned}
&-\mathrm{Q}_{n}\varphi_{0}(x)=e^{\frac{n}{n-1}\tau_{n}\varphi_{0}(x)}~\mathrm {in}~\mathbb{R}^{n},\\
&\varphi_{0}(x)\leq\varphi_{0}(0)=0,~x\in\mathbb{R}^{n}.
\end{aligned}
\right.\notag
\end{equation}
Furthermore, we obtain
\begin{equation}
\varphi_{0}(x)=-\frac{n-1}{\tau_{n}}\log(1+\kappa_{n}^{\frac{1}{n-1}}F^{\circ}(x)^{\frac{n}{n-1}})\notag
\end{equation}
and
\begin{equation}\label{E:5.31}
\int_{\mathbb{R}^{n}}^{}e^{\frac{n}{n-1}\tau_{n}\varphi_{0}(x)}dx=1.
\end{equation}

\textbf{Step 2.} The asymptotic behavior away from  $ x_{0 }$.

We have $ u_{\beta,M}=\min\{\frac{d_{\beta}}{M},u_{\beta}\} $ as recalled. For any fixed $R>0$, we have
\begin{align}\label{E:5.32}
\Vert u_{\beta,M}\Vert_{F(\Omega)}^{n}\geq&\frac{1}{\lambda_{\beta}}\int_{\mathcal{W}_{Rr_{\beta}^{\frac{1}{1-\beta}}}(x_{\beta})}^{}\frac{u_{\beta,M}u_{\beta}^{\frac{1}{n-1}}e^{\tau_{n}(1-\beta) u_{\beta}^{\frac{n}{n-1}}}}{F^{\circ}(x)^{n\beta}}dx\notag\\
=&\frac{1}{M}(1+o_{\beta}(1))\int_{\mathcal{W}_{R(0)}}^{}\frac{e^{\tau_{n}(1-\beta)(u_{\beta}^{\frac{n}{n-1}}(x_{\beta}+r_{\beta}^{\frac{1}{1-\beta}}y)-d_{\beta}^{\frac{n}{n-1}})}}{F^{\circ}\bigg(y+\frac{x_{\beta}}{r_{\beta}^{\frac{1}{1-\beta}}}\bigg)^{n\beta}}dy\notag\\
\to&\frac{1}{M}\int_{\mathcal{W}_{R}(0)}^{}e^{\frac{n}{n-1}\tau_{n}\varphi_{0}(y)}dy~ \mathrm {as}~\beta\to0
\end{align}
and
\begin{align}\label{E:5.33}
\Vert (u_{\beta}-\frac{d_{\beta}}{M})^{+}\Vert_{F(\Omega)}^{n}&\geq\frac{1}{\lambda_{\beta}}\int_{\mathcal{W}_{Rr_{\beta}^{\frac{1}{1-\beta}}}(x_{\beta})}^{}\frac{(u_{\beta}-\frac{d_{\beta}}{M})^{+}u_{\beta}^{\frac{1}{n-1}}e^{\tau_{n}(1-\beta) u_{\beta}^{\frac{n}{n-1}}}}{F^{\circ}(x)^{n\beta}}dx\notag\\
&=(1-\frac{1}{M})(1+o_{\beta}(1))\int_{\mathcal{W}_{R(0)}}^{}\frac{e^{\tau_{n}(1-\beta)(u_{\beta}^{\frac{n}{n-1}}(x_{\beta}+r_{\beta}^{\frac{1}{1-\beta}}y)-d_{\beta}^{\frac{n}{n-1}})}}{F^{\circ}\bigg(y+\frac{x_{\beta}}{r_{\beta}^{\frac{1}{1-\beta}}}\bigg)^{n\beta}}dy\notag\\
&\to(1-\frac{1}{M}) \int_{\mathcal{W}_{R}(0)}^{}e^{\frac{n}{n-1}\tau_{n}\varphi_{0}(y)}dy~\mathrm {as}~\beta\to0.
\end{align}
Combining \eqref{E:5.7}, \eqref{E:5.31}, \eqref{E:5.32} with \eqref{E:5.33}, and passing to the limits $\beta\rightarrow0 $, $ R\rightarrow+\infty $, we once again arrive at \eqref{E:5.8}.

Analogous to Lemma \ref{le:3.3}, we have the following result:
\begin{lem}\label{le:3.8}
There holds
\begin{equation}
\lim_{\beta\rightarrow0}\int_{\Omega}^{} \frac{e^{\tau_{n}(1-\beta) u_{\beta}^{\frac{n}{n-1}}}}{F^{\circ}(x)^{n\beta}}dx\leq\left|\Omega\right |+\lim_{R \to+\infty} \limsup\limits_{\beta\rightarrow0}\int_{\mathcal{W}_{Rr_{\beta}^{\frac{1}{1-\beta}}}(x_{\beta})}^{}\frac{e^{\tau_{n}(1-\beta)u_{\beta}^{\frac{n}{n-1}}}}{F^{\circ}(x)^{n\beta}}dx.\notag
\end{equation}
\end{lem}
\begin{proof}
For any $R>0$, elementary calculations give
\begin{align}
\int_{\mathcal{W}_{Rr_{\beta}^{\frac{1}{1-\beta}}}(x_{\beta})}^{}\frac{e^{\tau_{n}(1-\beta)u_{\beta}^{\frac{n}{n-1}}}}{F^{\circ}(x)^{n\beta}}dx&=\int_{\mathcal{W}_{R}(0)}^{}\frac{r_{\beta}^{n}e^{\tau_{n}(1-\beta)u_{\beta}^{\frac{n}{n-1}}(x_{\beta}+r_{\beta}^{\frac{1}{1-\beta}}y)}}{F^{\circ}\bigg(y+\frac{x_{\beta}}{r_{\beta}^{\frac{1}{1-\beta}}}\bigg)^{n\beta}}dy\notag\\
&=\frac{\lambda_{\beta}}{d_{\beta}^{\frac{n}{n-1}}}\int_{\mathcal{W}_{R}(0)}^{}\frac{e^{\tau_{n}(1-\beta)(u_{\beta}^{\frac{n}{n-1}}(x_{\beta}+r_{\beta}^{\frac{1}{1-\beta}}y)-d_{\beta}^{\frac{n}{n-1}})}}{F^{\circ}\bigg(y+\frac{x_{\beta}}{r_{\beta}^{\frac{1}{1-\beta}}}\bigg)^{n\beta}}dy.\notag\\
&\to\frac{\lambda_{\beta}}{d_{\beta}^{\frac{n}{n-1}}}\int_{\mathcal{W}_{R}(0)}^{}e^{\frac{n}{n-1}\tau_{n}\varphi_{0}(y)}dy~\mathrm {as}~\beta\to0.\notag
\end{align}
Hence
\begin{align}\label{E:5.34}
\lim\limits_{R\rightarrow+\infty} \limsup\limits_{\beta\rightarrow0}\int_{\mathcal{W}_{Rr_{\beta}^{\frac{1}{1-\beta}}}(x_{\beta})}^{}\frac{e^{\tau_{n}(1-\beta)u_{\beta}^{\frac{n}{n-1}}}}{F^{\circ}(x)^{n\beta}}dx=\limsup\limits_{\beta\rightarrow0}\frac{\lambda_{\beta}}{d_{\beta}^{\frac{n}{n-1}}}.\notag
\end{align}
This together with \eqref{E:5.12} implies that Lemma \ref{le:3.8} holds.
\end{proof}

\textbf{Step 3.}
We proceed to estimate the upper bound of $\mathfrak{T}$.
\begin{lem}\label{le:3.9}
There holds
\begin{equation}
\mathfrak{T}\leq\vert\Omega\vert+\kappa_{n}e^{\tau_{n}C_{0}+{\sum_{j=1}^{n-1}}\frac{1}{j}}.\notag
\end{equation}
\end{lem}
\begin{proof}
As before, set $ u_{\beta}(\rho)=\sup_{\partial \mathcal{W}_{\rho}}u_{\beta} $, $\tilde{u}_{\beta}=(u_{\beta}-u_{\beta}(\rho))^{+}$ and $s_{\beta}=\Vert \tilde{u}_{\beta}\Vert_{F(\mathcal{W}_{\rho})}$. Due to $r_{\beta}\to0$ as $\beta\to0$, we have $\mathcal{W}_{Rr_{\beta}^{1/(1-\beta)}}(x_{\beta})\subset \mathcal{W}_{\rho}$ for  $ \beta$ sufficently small. It is not difficult to see that 
\begin{align}\label{E:5.35}
\int_{\mathcal{W}_{Rr_{\beta}^{\frac{1}{1-\beta}}}(x_{\beta})}^{}\frac{e^{\tau_{n}(1-\beta)(\tilde{u}_{\beta}/s_{\beta})^{\frac{n}{n-1}}}}{F^{\circ}(x)^{n\beta}}dx
&\leq\int_{\mathcal{W}_{\rho}}^{}\frac{e^{\tau_{n}(1-\beta)(\tilde{u}_{\beta}/s_{\beta})^{\frac{n}{n-1}}}-1}{F^{\circ}(x)^{n\beta}}dx+o_{\beta}(1).
\end{align}
We conclude from Lemma \ref{le:3.6} and \eqref{E:5.35} that
\begin{align}
&\int_{\mathcal{W}_{Rr_{\beta}^{\frac{1}{1-\beta}}}(x_{\beta})}^{}\frac{e^{\tau_{n}(1-\beta)u_{\beta}^{\frac{n}{n-1}}}}{F^{\circ}(x)^{n\beta}}dx\notag\\
\leq&\int_{\mathcal{W}_{Rr_{\beta}^{\frac{1}{1-\beta}}}(x_{\beta})}^{}\frac{e^{\tau_{n}(1-\beta)((\tilde{u}_{\beta}/s_{\beta})^{\frac{n}{n-1}}-\frac{n}{\tau_{n}}\log\rho+C_{0}+o_{\beta,\rho}(1))}}{F^{\circ}(x)^{n\beta}}dx\notag\\
\leq&\rho^{-n(1-\beta)}e^{\tau_{n}(1-\beta) C_{0}+o_{\beta,\rho}(1)}\bigg(\int_{\mathcal{W}_{\rho}}^{}\frac{e^{\tau_{n}(1-\beta)(\tilde{u}_{\beta}/s_{\beta})^{\frac{n}{n-1}}}-1}{F^{\circ}(x)^{n\beta}}dx+o_{\beta}(1)\bigg).\notag\\
\leq&\rho^{n\beta}\kappa_{n}e^{\tau_{n}(1-\beta) C_{0}+{\sum_{j=1}^{n-1}}\frac{1}{j}}+o_{\beta,\rho}(1).\notag
\end{align}
Therefore,
\begin{equation}\label{E:3.401}
\lim_{R\to+\infty}\limsup\limits_{\beta\rightarrow0}\int_{\mathcal{W}_{Rr_{\beta}^{\frac{1}{1-\beta}}}(x_{\beta})}^{}\frac{e^{\tau_{n}(1-\beta)u_{\beta}^{\frac{n}{n-1}}}}{F^{\circ}(x)^{n\beta}}dx\leq\kappa_{n}e^{\tau_{n}C_{0}+{\sum_{j=1}^{n-1}}\frac{1}{j}}.
\end{equation}
From \eqref{E:3.3}, \eqref{E:3.401} and Lemma \ref{le:3.8}, we arrive at Lemma \ref{le:3.9}.
\end{proof}
\vspace{0.2em}
\subsection{Case (ii): $x_{0}\in\Omega$ and $x_{0}\neq0$}
In this case, we will again repeat the three steps from Case (i).

\textbf{Step 1.} The asymptotic behavior of $u_{\beta}$ near $x_{0}$.

We set 
\begin{equation}\label{E:5.361}
\chi_{\beta}(x)=d_{\beta}^{-1}u_{\beta}(x_{\beta}+r_{\beta}x) ~\mathrm {and}~ \psi_{\beta}(x)=d_{\beta}^{\frac{1}{n-1}}(u_{\beta}(x_{\beta}+r_{\beta}x)-d_{\beta})
\end{equation}
for $x\in \Omega_{\beta,3}=\{x\in\mathbb{R}^{n}:x_{\beta}+r_{\beta}x\in\Omega\}$, where $\Omega_{\beta,3}\to\mathbb{R}^{n}$ as $\beta\to0$. A direct computation yields
\begin{equation}\label{E:5.36}
-\mathrm{Q}_{n}\chi_{\beta}(x)=d_{\beta}^{-n}F^{\circ}(x_{\beta}+r_{\beta}x)^{-n\beta}\chi_{\beta}^{\frac{1}{n-1}}e^{\tau_{n}(1-\beta)(u_{\beta}^{\frac{n}{n-1}}(x_{\beta}+r_{\beta}x)-d_{\beta}^{\frac{n}{n-1}})}~\mathrm{in}~ \Omega_{\beta,3},
\end{equation}
and
\begin{equation}\label{E:5.37}
-\mathrm{Q}_{n}\psi_{\beta}(x)=F^{\circ}(x_{\beta}+r_{\beta}x)^{-n\beta}\chi_{\beta}^{\frac{1}{n-1}}e^{\tau_{n}(1-\beta)(u_{\beta}^{\frac{n}{n-1}}(x_{\beta}+r_{\beta}x)-d_{\beta}^{\frac{n}{n-1}})}~ \mathrm{in}~ \Omega_{\beta,3}. 
\end{equation}
Since $r_{\beta}\to0$ and $x_{\beta}\not=0$, we conclude that $F^{\circ}(x_{\beta}+r_{\beta}x)^{-n\beta}$ is bounded in $L^{p}(\mathcal{W}_{R})$ for some $p>1$ and fixed $R>0$. Therefore, we have $\mathrm{Q}_{n}\chi_{\beta}(x)$ and $\mathrm{Q}_{n}\psi_{\beta}(x)$ is bounded in $L^{p}(\mathcal{W}_{R})$ for some $p>1$. By employing the Regularity Theory to \eqref{E:5.36} and \eqref{E:5.37}, we obtain 
\begin{equation}
\left\{  
\begin{aligned}
&\chi_{\beta}(x)\rightarrow\chi_{0}(x)~ \mathrm {in}~ C_{\mathrm {loc}}^{1}(\mathbb{R}^{n})~\mathrm {as}~\beta\to0,\\
&\psi_{\beta}(x)\rightarrow\psi_{0}(x)~ \mathrm {in}~ C_{\mathrm {loc}}^{1}(\mathbb{R}^{n})~\mathrm {as}~\beta\to0,&\notag
	\end{aligned}
	\right.
\end{equation}
where  $\chi_{0}(x)=\chi_{0}(0)\equiv1 $, $\psi_{0}(0)=0=\mathrm{max}_{\mathbb{R}^{n}}\psi_{0}(x)$ and $\psi_{0}(x)$ solves $$-\mathrm{Q}_{n}\psi_{0}(x)=e^{\frac{n}{n-1}\tau_{n}\psi_{0}(x)}$$ in the distributional sense. 
Furthermore, we have  
\begin{equation}\label{E:5.391}
\psi_{0}(x)=-\frac{n-1}{\tau_{n}}\log(1+\kappa_{n}^{\frac{1}{n-1}}F^{\circ}(x)^{\frac{n}{n-1}}),
\end{equation}
and
\begin{equation}\label{E:5.38}
\int_{\mathbb{R}^{n}}^{}e^{\frac{n}{n-1}\tau_{n}\psi_{0}(x)}dx=1.
\end{equation}

\textbf{Step 2.} 
The asymptotic behavior away from  $ x_{0 }$.

Direct calculation shows that 
\begin{align}\label{E:5.39}
\Vert u_{\beta,M}\Vert_{F(\Omega)}^{n}\geq&\frac{1}{\lambda_{\beta}}\int_{\mathcal{W}_{Rr_{\beta}}(x_{\beta})}^{}\frac{u_{\beta,M}u_{\beta}^{\frac{1}{n-1}}e^{\tau_{n}(1-\beta) u_{\beta}^{\frac{n}{n-1}}}}{F^{\circ}(x)^{n\beta}}dx\notag\\
=&\frac{1}{M}(1+o_{\beta}(1))\int_{\mathcal{W}_{R(0)}}^{}\frac{e^{\tau_{n}(1-\beta)(u_{\beta}^{\frac{n}{n-1}}(x_{\beta}+r_{\beta}y)-d_{\beta}^{\frac{n}{n-1}})}}{F^{\circ}(x_{\beta}+r_{\beta}y)^{n\beta}}dy\notag\\
\to&\frac{1}{M}\int_{\mathcal{W}_{R(0)}}^{}e^{\frac{n}{n-1}\tau_{n}\psi_{0}(y)}dy~ \mathrm {as}~\beta\to0
\end{align}
and
\begin{align}\label{E:5.40}
\Vert(u_{\beta}-\frac{d_{\beta}}{M})^{+}\Vert_{F(\Omega)}^{n}&\geq\frac{1}{\lambda_{\beta}}\int_{\mathcal{W}_{Rr_{\beta}}(x_{\beta})}^{}\frac{(u_{\beta}-\frac{d_{\beta}}{M})^{+}u_{\beta}^{\frac{1}{n-1}}e^{\tau_{n}(1-\beta) u_{\beta}^{\frac{n}{n-1}}}}{F^{\circ}(x)^{n\beta}}dx\notag\\
&=(1-\frac{1}{M})(1+o_{\beta}(1))\int_{\mathcal{W}_{R}}^{}\frac{e^{\tau_{n}(1-\beta)(u_{\beta}^{\frac{1}{n-1}}(x_{\beta}+r_{\beta}y)-d_{\beta}^{\frac{1}{n-1}})}}{F^{\circ}(x_{\beta}+r_{\beta}y)^{n\beta}}dy\notag\\
&\to(1-\frac{1}{M}) \int_{\mathcal{W}_{R(0)}}^{}e^{\frac{n}{n-1}\tau_{n}\psi_{0}(y)}dy~\mathrm {as}~\beta\to0.
\end{align}
In view of \eqref{E:5.7}, \eqref{E:5.38}, \eqref{E:5.39} and \eqref{E:5.40}, we again obtain \eqref{E:5.8}.

For any $R>0$, one gets
\begin{align}
\int_{\mathcal{W}_{Rr_{\beta}}(x_{\beta})}^{}\frac{e^{\tau_{n}(1-\beta)u_{\beta}^{\frac{n}{n-1}}}}{F^{\circ}(x)^{n\beta}}dx&=\int_{\mathcal{W}_{R(0)}}^{}\frac{r_{\beta}^{n}e^{\tau_{n}(1-\beta)u_{\beta}^{\frac{n}{n-1}}(x_{\beta}+r_{\beta}y)}}{F^{\circ}(x_{\beta}+r_{\beta}y)^{n\beta}}dy\notag\\
&=\frac{\lambda_{\beta}}{d_{\beta}^{\frac{n}{n-1}}}\int_{\mathcal{W}_{R(0)}}^{}\frac{e^{\tau_{n}(1-\beta)(u_{\beta}^{\frac{n}{n-1}}(x_{\beta}+r_{\beta}y)-d_{\beta}^{\frac{n}{n-1}})}}{F^{\circ}(x_{\beta}+r_{\beta}y)^{n\beta}}dy\notag\\
&\to\frac{\lambda_{\beta}}{d_{\beta}^{\frac{n}{n-1}}}\int_{\mathcal{W}_{R(0)}}^{}e^{\frac{n}{n-1}\tau_{n}\psi_{0}(y)}dy~\mathrm {as}~\beta\to0.\notag
\end{align}
Then we have
\begin{align}\label{E:5.41}
\lim\limits_{R\rightarrow+\infty}\limsup\limits_{\beta\rightarrow0}\int_{\mathcal{W}_{Rr_{\beta}}(x_{\beta})}^{}\frac{e^{\tau_{n}(1-\beta)u_{\beta}^{\frac{n}{n-1}}}}{F^{\circ}(x)^{n\beta}}dx=&\limsup\limits_{\beta\rightarrow0}\frac{\lambda_{\beta}}{d_{\beta}^{\frac{n}{n-1}}}.
\end{align}
It follows from \eqref{E:5.12} and \eqref{E:5.41} that
\begin{equation}\label{E:5.42}
\lim_{\beta\rightarrow0}\int_{\Omega}^{} \frac{e^{\tau_{n}(1-\beta) u_{\beta}^{\frac{n}{n-1}}}}{F^{\circ}(x)^{n\beta}}dx\leq\left|\Omega\right |+\lim_{R \to+\infty} \limsup\limits_{\beta\rightarrow0}\int_{\mathcal{W}_{Rr_{\beta}}(x_{\beta})}^{}\frac{e^{\tau_{n}(1-\beta)u_{\beta}^{\frac{n}{n-1}}}}{F^{\circ}(x)^{n\beta}}dx.
\end{equation}

Furthermore, we have $ d_{\beta}^{\frac{1}{n-1}}u_{\beta}\rightharpoonup G_{x_{0}}$ weakly in $W_{0}^{1,q}(\Omega)$ for any $ 1<q<n $ and $ d_{\beta}^{\frac{1}{n-1}}u_{\beta}\rightarrow G_{x_{0}}$ in $ C_{\mathrm{loc}}^{1}(\overline{\Omega}\setminus\{x_{0}\}) $, where $G_{x_{0}}$ satisfies $ -\mathrm{Q}_{n}G_{x_{0}}=\delta_{x_{0}}$ in a distributional sense. Moreover, $G_{x_{0}}$ takes the form
\begin{equation}\label{E:5.451}
G_{x_{0}}(x)=-\frac{n}{\tau_{n}}\log F^{\circ}(x-x_{0})+C_{x_{0}}+ D(x),
\end{equation} 
where $C_{x_{0}}$ is a constant, $ D(x)\in C^{0}(\overline{\Omega})\cap C_{\mathrm{loc}}^{1}(\overline{\Omega}\setminus\{x_{0}\})$ satisfies $D(x_{0})=0 $ and $\lim\limits_{x\rightarrow x_{0}}F^{\circ}(x-x_{0})\nabla D(x)=0$. 

\textbf{Step 3.}
The upper bound for $\mathfrak{T}$.
\begin{lem}\label{le:3.10}
There holds
\begin{equation}
\mathfrak{T}\leq\vert\Omega\vert+\kappa_{n}e^{\tau_{n}C_{x_{0}}+{\sum_{j=1}^{n-1}}\frac{1}{j}}.\notag
\end{equation}
\end{lem}
\begin{proof}
We write $ u_{\beta}(\rho)=\sup_{\partial \mathcal{W}_{\rho}(x_{0})}u_{\beta} $, $\tilde{u}_{\beta}=(u_{\beta}-u_{\beta}(\rho))^{+}$ and $s_{\beta}=\Vert \tilde{u}_{\beta}\Vert_{F(\mathcal{W}_{\rho}(x_{0}))}$. Note that $\mathcal{W}_{Rr_{\beta}(x_{0})}\subset \mathcal{W}_{\rho}(x_{0})$ for $ \beta$ sufficently small. By direct computation, we find
\begin{align}
\int_{\mathcal{W}_{Rr_{\beta}}(x_{\beta})}^{}\frac{e^{\tau_{n}(1-\beta)u_{\beta}^{\frac{n}{n-1}}}}{F^{\circ}(x)^{n\beta}}dx\leq\rho^{-n\beta}\kappa_{n}e^{\tau_{n}(1-\beta)C_{x_{0}}+{\sum_{j=1}^{n-1}}\frac{1}{j}}+o_{\beta,\rho}(1).\notag
\end{align}
Letting $ \beta\rightarrow0 $ and $ R\rightarrow+\infty $, one has
\begin{equation}
\lim_{R\to+\infty}\limsup\limits_{\beta\rightarrow0}\int_{\mathcal{W}_{Rr_{\beta}}(x_{\beta})}^{}\frac{e^{\tau_{n}(1-\beta)u_{\beta}^{\frac{n}{n-1}}}}{F^{\circ}(x)^{n}}dx\leq\kappa_{n}e^{\tau_{n}C_{x_{0}}+{\sum_{j=1}^{n-1}}\frac{1}{j}}.\notag 
\end{equation}
This together with \eqref{E:3.3} and \eqref{E:5.42} leads to  Lemma \ref{le:3.10}.
\end{proof}

\subsection{Case (iii): $x_{0}\in\partial\Omega$}
In this case, we set $\chi_{\beta}$ and $\psi_{\beta}$ be defined by \eqref{E:5.361}. According to Lemma 4.7 in \cite{ref22}, we have the following result:  
\begin{lem}\label{le:3.11}
As $\beta\to0$, we have $\frac{dist(x_{\beta},\partial\Omega)}{r_{\beta}}\to+\infty$, $ \chi_{\beta}(x)\rightarrow1 $ and $ \psi_{\beta}(x)\rightarrow\psi_{0}(x) $ in $ C_{\mathrm {loc}}^{1}(\mathbb{R}^{n}) $. Consequently, \eqref{E:5.391} and \eqref{E:5.38} remain valid. In addition, $d_{\beta}^{\frac{1}{n-1}}u_{\beta}\rightharpoonup\tilde{G}$ weakly in $W_{0}^{1,n}(\Omega)$ and $d_{\beta}^{\frac{1}{n-1}}u_{\beta}\to\tilde{G}$ in $C^{1}(\overline{\Omega}\setminus\left \{ x_{0} \right \})$, where $\tilde{G}=0$ solves
\begin{equation}
\left\{  
\begin{aligned}
&-\mathrm{Q}_{n}\tilde{G}=0,~\mathrm {in}~\Omega,\\
&\tilde{G}=0,~ \mathrm {in}~\partial\Omega,&\notag
\end{aligned}
\right.
\end{equation}
uniquely.
\end{lem}

Writing $u_{\beta}(\rho)=\sup_{\partial \mathcal{W}_{\rho}(x_{0})}u_{\beta} $ and $\tilde{u}_{\beta}=(u_{\beta}-u_{\beta}(\rho))^{+}$. From Lemma \ref{le:3.11}, it follows that
\begin{equation}\label{E:5.45}
u_{\beta}(\rho)=\frac{1}{d_{\beta}^{\frac{1}{n-1}}}\bigg(\sup_{\partial \mathcal{W}_{\rho}(x_{0})}\tilde{G}+o_{\beta}(1)\bigg)=\frac{o_{\beta}(1)}{d_{\beta}^{\frac{1}{n-1}}}
\end{equation}
and 
\begin{equation}
\Vert u_{\beta}\Vert_{F(\Omega\setminus \mathcal{W}_{\rho}(x_{0}))}^{n}=\frac{1}{d_{\beta}^{\frac{n}{n-1}}}\int_{\Omega\setminus \mathcal{W}_{\rho}(x_{0})}^{}\left(F^{n}(\nabla \tilde{G})+o_{\beta}(1)\right)dx=\frac{o_{\beta,\rho}(1)}{d_{\beta}^{\frac{n}{n-1}}}.\notag
\end{equation}
Therefore, we have
\begin{equation}\label{E:5.46}
s_{\beta}=\Vert \tilde{u}_{\beta}\Vert_{F(\mathcal{W}_{\rho}(x_{0}))}=\bigg(1-\frac{o_{\beta,\rho}(1)}{d_{\beta}^{\frac{n}{n-1}}}\bigg)^{1/n}.
\end{equation}
Then using \eqref{E:5.45} and \eqref{E:5.46}, we get
\begin{align}\label{E:5.47}
\tau_{n}(1-\beta) u_{\beta}^{\frac{n}{n-1}}\notag\leq&\tau_{n}(1-\beta)\bigg(\tilde{u}_{\beta}^{\frac{n}{n-1}}+\frac{n}{n-1}\tilde{u}_{\beta}^{\frac{1}{n-1}}u_{\beta}(\delta)+o_{\beta}(1)\bigg) \notag\\
\leq&\tau_{n}(1-\beta)\bigg(\frac{\tilde{u}_{\beta}}{s_{\beta}}\bigg)^{\frac{n}{n-1}}+o_{\beta}(1) 
\end{align}
for $x\in \mathcal{W}_{Rr_{\beta}(x_{\beta})}$. It follows from \eqref{E:5.47} that
\begin{align}
\int_{\mathcal{W}_{Rr_{\beta}}(x_{\beta})}^{}\frac{e^{\tau_{n}(1-\beta)u_{\beta}^{\frac{n}{n-1}}}}{F^{\circ}(x)^{n\beta}}dx
\leq&\int_{\mathcal{W}_{Rr_{\beta}}(x_{\beta})}^{}\frac{e^{\tau_{n}(1-\beta)(\tilde{u}_{\beta}/s_{\beta})^{\frac{n}{n-1}}}}{F^{\circ}(x)^{n\beta}}dx+o_{\beta,R}(1)\notag\\
	\leq&\int_{\mathcal{W}_{\rho}(x_{0})}^{}\frac{e^{\tau_{n}(1-\beta)(\tilde{u}_{\beta}/s_{\beta})^{\frac{n}{n-1}}}-1}{F^{\circ}(x)^{n\beta}}dx+o_{\beta,R}(1)\notag\\
	\leq&\rho^{n}\kappa_{n}e^{{\sum_{j=1}^{n-1}}\frac{1}{j}}+o_{\beta,R}(1).\notag
\end{align}
Thus
\begin{equation}\label{E:5.48}
\lim_{R\to+\infty}\limsup\limits_{\beta\rightarrow0}\int_{\mathcal{W}_{Rr_{\beta}}(x_{\beta})}^{}\frac{e^{\tau_{n}(1-\beta)u_{\beta}^{\frac{n}{n-1}}}}{F^{\circ}(x)^{n\beta}}dx\leq\rho^{n}\kappa_{n}e^{{\sum_{j=1}^{n-1}}\frac{1}{j}}.
\end{equation}
Taking the limit $\rho\to0$ in the combination of \eqref{E:3.3}, \eqref{E:5.42} and \eqref{E:5.48}, we obtain
\begin{equation}
	\mathfrak{T}\leq\vert\Omega\vert.\notag
\end{equation} 
This contradicts
$\mathfrak{T}>\vert \Omega\vert$, hence case (iii) cannot occur.

\textbf{Conclusion:} Up to now, we have completed the blow-up analysis and concluded that the blow-up point $x_{0}$ can only lie in the interior of $\Omega$. That is, cases (i) and (ii) hold. Moreover, in case (i), Lemmas \ref{le:3.7} and \ref{le:3.9} establish upper bounds for $\mathfrak{T}$, while Lemma \ref{le:3.10} provides the corresponding upper bound for case (ii). More precisely, we have
\begin{equation}\label{E:3.544}
\mathfrak{T}\leq\vert\Omega\vert+\kappa_{n}e^{\tau_{n}C_{x_{0}}+{\sum_{j=1}^{n-1}}\frac{1}{j}}.
\end{equation}

\section{Proof of Theorem \ref{th:1.1}}
In this section, we will find a test function $\upsilon_{\beta}(x)\in W_{0}^{1,n}(\Omega)$ satisfying $\Vert\upsilon_{\beta}(x)\Vert_{F(\Omega)}=1$ and 
\begin{equation}\label{E:6.1}
\int_{\Omega}^{} e^{\tau_{n}\vert\upsilon_{\beta}(x)\vert^{\frac{n}{n-1}}}dx>\vert\Omega\vert+\kappa_{n}e^{\tau_{n}C_{x_{0}}+{\sum_{j=1}^{n-1}}\frac{1}{j}},
\end{equation}
which implies 
\begin{equation}
\mathfrak{T}>\vert\Omega\vert+\kappa_{n}e^{\tau_{n}C_{x_{0}}+{\sum_{j=1}^{n-1}}\frac{1}{j}}.\notag
\end{equation}
This contradicts \eqref{E:3.544}, thus we eliminate the blow-up phenomenon. Since $u_{\beta}$ is uniformly bounded in $\Omega$, applying elliptic estimates to \eqref{E:3.1}, we can obtain $u_{\beta}\to u_{0}$ in $C^{1}(\overline{\Omega})$ as $\beta\to0$. The Lebesgue Dominated Convergence Theorem states that
\begin{equation}
	\int_{\Omega}^{} e^{\tau_{n} u_{0}^{\frac{n}{n-1}}}dx= \lim\limits_{\beta\rightarrow0}\int_{\Omega}^{} \frac{e^{\tau_{n}(1-\beta) u_{\beta}^{\frac{n}{n-1}}}}{F^{\circ}(x)^{n\beta}}dx=\mathfrak{T}.\notag
\end{equation}
In this way, we complete the proof of Theorem \ref{th:1.1}.
 
Now we are going to construct the test functions $\upsilon_{\beta}(x)\in W_{0}^{1,n}(\Omega)$ satisfying \eqref{E:6.1}. Let $R=-\log\beta$ and $G_{x_{0}}$ given by \eqref{E:5.451}. We define
\begin{equation}
\upsilon_{\beta}(x)=\left\{
\begin{aligned}
&d+\frac{1}{d^{\frac{1}{n-1}}}\bigg(-\frac{n-1}{\tau_{n}}\log\bigg(1+\kappa_{n}^{\frac{1}{n-1}}\bigg(\frac{F^{\circ}(x-x_{0})}{\beta}\bigg)^{\frac{n}{n-1}}\bigg)+m\bigg),~x\in \mathcal{W}_{R\beta}(x_{0}),\\
&\frac{G_{x_{0}}}{d^{\frac{1}{n-1}}},~ x\in\Omega\setminus \mathcal{W}_{R\beta}(x_{0}),&
	\end{aligned}
	\right.\notag
\end{equation}
where $d$ and $m$ are constants depending on $\beta$. To ensure $\upsilon_{\beta}(x)\in W_{0}^{1,n}(\Omega)$, we set
\begin{equation}
d+\frac{1}{d^{\frac{1}{n-1}}}\bigg(-\frac{n-1}{\tau_{n}}\log\bigg(1+\kappa_{n}^{\frac{1}{n-1}}R^{\frac{n}{n-1}}\bigg)+m\bigg)=\frac{G_{x_{0}}(R\beta)}{d^{\frac{1}{n-1}}}.\notag
\end{equation}
For any $R>0$, since $R\beta\to0$ as $\beta\to0$, it follows that 
\begin{equation}
G_{x_{0}}(R\beta)=-\frac{n}{\tau_{n}}\log R-\frac{n}{\tau_{n}}\log\beta+C_{x_{0}}+o_{\beta}(1).\notag
\end{equation}
Calculations show that
\begin{equation}\label{E:6.2}
d^{\frac{n}{n-1}}=-\frac{n}{\tau_{n}}\log\beta+C_{x_{0}}+\frac{1}{\tau_{n}}\log\kappa_{n}-m+O\bigg(\frac{1}{R^{\frac{n}{n-1}}}\bigg).
\end{equation}

A straight calculation gives
\begin{align}
\Vert\upsilon_{\beta}(x)\Vert_{F(\mathcal{W}_{R\beta}(x_{0}))}^{n}=&\frac{1}{n^{\frac{n}{n-1}}\beta^{\frac{n^{2}}{n-1}}d^{\frac{n}{n-1}}}\int_{\mathcal{W}_{R\beta}(x_{0})}^{}\frac{F^{\circ}(x-x_{0})^{\frac{n}{n-1}}}{(1+\kappa_{n}^{\frac{1}{n-1}}(\frac{F^{\circ}(x-x_{0})}{\beta})^{\frac{n}{n-1}})^{n}}dx\notag\\
=&\frac{1}{n^{\frac{n}{n-1}}\beta^{\frac{n^{2}}{n-1}}d^{\frac{n}{n-1}}}\int_{0}^{R\beta}\frac{n\kappa_{n}r^{\frac{n^{2}-n+1}{n-1}}}{(1+\kappa_{n}^{\frac{1}{n-1}}(\frac{r}{\beta})^{\frac{n}{n-1}})^{n}}dr\notag\\
=&\frac{n-1}{\tau_{n}d^{\frac{n}{n-1}}}\int_{1}^{1+\kappa_{n}^{\frac{1}{n-1}}R^{\frac{n}{n-1}}}\frac{(z-1)^{n-1}}{z^{n}}dz\notag\\
=&\frac{1}{\tau_{n}d^{\frac{n}{n-1}}}\bigg(\log\kappa_{n}+n\log R- (n-1){\textstyle \sum_{j=1}^{n-1}}\frac{1}{j}+O\bigg(\frac{1}{R^{\frac{n}{n-1}}}\bigg)\bigg),\notag
\end{align}
where $r=F^{\circ}(x-x_{0})$ and $z=1+\kappa_{n}^{\frac{1}{n-1}}(\frac{r}{\beta})^{\frac{n}{n-1}}$.
Meanwhile, one gets
\begin{align}
\Vert\upsilon_{\beta}(x)\Vert_{F(\Omega\setminus \mathcal{W}_{R\beta}(x_{0}))}^{n}=&\frac{1}{d^{\frac{n}{n-1}}}\bigg(\int_{\Omega\setminus \mathcal{W}_{R\beta}(x_{0})}^{}F^{n}(\nabla G_{x_{0}})dx+o_{\beta,R}(1)\bigg)\notag\\
=&\frac{1}{d^{\frac{n}{n-1}}}\bigg(-\int_{\partial \mathcal{W}_{R\beta}(x_{0})}^{}G_{x_{0}}F^{^{n-1}}(\nabla G_{x_{0}})\left \langle F_{\xi}(\nabla G_{x_{0}}),\vec{n} \right\rangle ds+o_{\beta,R}(1)\notag\\&-\int_{\Omega\setminus \mathcal{W}_{R\beta}(x_{0})}^{}G_{x_{0}}\mathrm{Q}_{n}G_{x_{0}}dx\bigg)\notag\\
=&\frac{1}{\tau_{n}d^{\frac{n}{n-1}}}(-n\log R-n\log\beta+\tau_{n}C_{x_{0}}+o_{\beta,R}(1)),\notag
\end{align}
where $ \vec{n} $ denotes the unit outward normal to $\partial \mathcal{W}_{R\beta}(x_{0})$.
To ensure $\Vert \upsilon_{\beta}(x)\Vert_{F(\Omega)}=1$, we set
\begin{equation}
\Vert\upsilon_{\beta}(x)\Vert_{F(\mathcal{W}_{R\beta}(x_{0}))}^{n}+\Vert\upsilon_{\beta}(x)\Vert_{F(\Omega\setminus \mathcal{W}_{R\beta}(x_{0}))}^{n}=1.\notag
\end{equation}
This indicates that
\begin{equation}\label{E:6.3}
d^{\frac{n}{n-1}}=-\frac{n-1}{\tau_{n}}{\textstyle \sum_{j=1}^{n-1}}\frac{1}{j}-\frac{n}{\tau_{n}}\log \beta+\frac{1}{\tau_{n}}\log\kappa_{n} +C_{x_{0}}+O\bigg(\frac{1}{R^{\frac{n}{n-1}}}\bigg)+o_{\beta,R}(1).
\end{equation}
Combining \eqref{E:6.2} and \eqref{E:6.3}, we obtain
\begin{equation}
m=\frac{n-1}{\tau_{n}}{\textstyle \sum_{j=1}^{n-1}}\frac{1}{j}+O\bigg(\frac{1}{R^{\frac{n}{n-1}}}\bigg)+o_{\beta,R}(1).\notag
\end{equation}
One can see that
\begin{align}
\tau_{n}\vert\upsilon_{\beta}(x)\vert^{\frac{n}{n-1}}\geq&-n\log\bigg(1+\kappa_{n}^{\frac{1}{n-1}}\bigg(\frac{F^{\circ}(x-x_{0})}{\beta}\bigg)^{\frac{n}{n-1}}\bigg)+{\textstyle \sum_{j=1}^{n-1}}\frac{1}{j}+\tau_{n}C_{x_{0}}+\log\kappa_{n}\notag\\
&-n\log\beta+O\bigg(\frac{1}{R^{\frac{n}{n-1}}}\bigg)+o_{\beta,R}(1)\notag
\end{align}
on $\mathcal{W}_{R\beta}(x_{0})$.
Then we get
\begin{align}\label{E:6.31}
\int_{\mathcal{W}_{R\beta}(x_{0})}^{} e^{\tau_{n}\vert\upsilon_{\beta}(x)\vert^{\frac{n}{n-1}}}dx\geq&e^{{ \sum_{j=1}^{n-1}}\frac{1}{j}+\tau_{n}C_{x_{0}}+\log\kappa_{n}-n\log\beta+O(\frac{1}{R^{\frac{n}{n-1}}})+o_{\beta,R}(1)}\notag\\
&\times\int_{\mathcal{W}_{R\beta}(x_{0})}^{}\frac{1}{(1+\kappa_{n}^{\frac{1}{n-1}}(\frac{F^{\circ}(x-x_{0})}{\beta})^{\frac{n}{n-1}})^{n}}dx\notag\\
\geq&\kappa_{n}e^{\tau_{n}C_{x_{0}}+{ \sum_{j=1}^{n-1}}\frac{1}{j}}+O\bigg(\frac{1}{R^{\frac{n}{n-1}}}\bigg)+o_{\beta,R}(1).
\end{align}
Meanwhile, we obtain that
\begin{align}
\int_{\Omega\setminus \mathcal{W}_{R\beta}(x_{0})}^{}e^{\tau_{n}\vert\upsilon_{\beta}(x)\vert^{\frac{n}{n-1}}}dx\geq&\int_{\Omega\setminus \mathcal{W}_{R\beta}(x_{0})}^{}(1+\tau_{n}\vert\upsilon_{\beta}(x)\vert^{\frac{n}{n-1}})dx\notag\\
=&\vert\Omega\vert+\frac{\tau_{n}}{d^{\frac{n}{(n-1)^{2}}}}\int_{\Omega\setminus \mathcal{W}_{R\beta}(x_{0})}^{}\vert G_{x_{0}}\vert^{\frac{n}{n-1}}dx,\notag
\end{align}
which together with \eqref{E:6.31} leads to \eqref{E:6.1}.

\bibliographystyle{plain}

\begin{thebibliography}{99}
\bibitem{ref11}Adimurthi, Sandeep K. A singular Moser-Trudinger embedding and its applications[J]. Nodea Nonlinear Differential Equations and Applications NoDEA, 2007, 5(13): 585-603.	
	
\bibitem{ref21} Alvino A, Ferone V, Trombetti G, Lions P L. Convex symmetrization and applications[J].
Ann. Inst. H. Poincar´e Anal. Non Lin´eaire, 1997, 14(2): 275-293.	
	
\bibitem{ref37} Bellettini G, Paolini M. Anisotroptic motion by mean curvature in the context of Finsler geometry[J]. Hokkaido Mathematical Journal, 1996, 25(3): 537-566.	
	
\bibitem{ref7}Carleson L, Chang S Y A. On the existence of an extremal function for an inequality of J. Moser[J]. Bulletin des Sciences Mathématiques, 1986, 110(2): 113-127.

\bibitem{ref12}Csat\'{o} G, Roy P. Extremal functions for the singular Moser-Trudinger inequality in 2 dimensions[J]. Calculus of Variations and Partial Differential Equations, 2015, 54(2): 2341-2366.	

\bibitem{ref14}Csat\'{o} G, Roy P, Nguyen V H. Extremals for the singular Moser-Trudinger inequality via n-harmonic transplantation[J]. Journal of Differential Equations, 2021, 270: 843-882.	

\bibitem{ref42}de Figueiredo D G, do \'{O} J M, dos Santos E M. Trudinger-Moser inequalities involving fast growth and weights with strong vanishing at zero[J]. Proceedings of the American Mathematical Society, 2016, 144(8): 3369-3380.	

\bibitem{ref29}Ding W Y, Jost J, Li J Y, Wang G F. The differential equation $ -\Delta u=8\pi-8\pi he^{u} $ on a compact Riemann Surface[J]. Asian Journal of Mathematics, 1997, 1(2): 230-248.

\bibitem{ref9}Flucher M. Extremal functions for Trudinger-Moser inequality in 2 dimensions[J]. Commentarii Mathematici Helvetici, 1992, 67(1): 471-497.	
	
\bibitem{ref28}Guo K, Liu Y. Sharp anisotropic singular Trudinger–Moser inequalities in the entire space[J]. Calculus of Variations and Partial Differential Equations, 2024, 63(4): 82.	
	
\bibitem{ref41} Heinonen J, Kilpelainen T, Martio O. Nonlinear Potential Theory of Degenerate Elliptic Equations. Oxford University Press, 1993.	

\bibitem{ref24}Li X M. Anisotropic singular Trudinger-Moser inequalities on the whole Euclidean space[J]. Discrete and Continuous Dynamical Systems, 2024, 44(2): 462-490.	

\bibitem{ref19}Li X M, Yang L, Su X F. Compactness of extremals for singular Trudinger-Moser inequalities one the whole Euclidean space[J]. Journal of Mathematical Inequalities, 2022, 16(3): 1179-1200.

\bibitem{ref34}Li X M, Yang Y Y. Extremal functions for singular Trudinger-Moser inequalities in the entire Euclidean space[J]. Journal
of Differential Equations, 2018, 264(8): 4901-4943.

\bibitem{ref30}Li Y X. Moser-Trudinger inequality on compact Riemannian manifolds of dimension two[J]. Partial Differential Equations, 2001, 14(2): 163-192.

\bibitem{ref31}Li Y X. Extremal functions for the Moser-Trudinger inequalities on compact Riemannian manifolds[J]. Science in China Series A: Mathematics, 2005, 48(5): 618-648.

\bibitem{ref38}Lieberman G M. Boundary regularity for solutions of degenerate elliptic equations[J]. Nonlinear Analysis: Theory, Methods Applications, 1988, 12(11): 1203-1219.

\bibitem{ref10}Lin K. Extremal functions for Moser's inequality[J]. Transactions of the American Mathematical Society, 1996, 348(7): 2663-2671.

\bibitem{ref27}Liu Y J. Concentration-compactness principle of singular Trudinger-Moser inequality involving N-Finsler-Laplacian operator[J]. International Journal of Mathematics, 2020 31(11) : 2050085

\bibitem{ref1}Lu G Z, Shen Y S, Xue J W, Zhu M C. Weighted anisotropic isoperimetric inequalities and existence of extremals for singular anisotropic Trudinger-Moser inequalities[J]. Advances in Mathematics, 2024, 458: 109949.

\bibitem{ref32}Lu G Z, Yang Y Y. The sharp constant and extremal functions for Moser-Trudinger inequalities involving $L^{p}$ norms[J]. Discrete and Continuous Dynamical Systems, 2009, 25(3): 963-979.	
	
\bibitem{ref18}Luo Q J, Li X M. The Compactness of Extremals for a Singular Hardy-Trudinger-Moser Inequality[J]. Journal of Partial Differential Equations, 2024, 37(3): 235-250.		
	
\bibitem{ref6}Moser J. A sharp form of an inequality by $ N $. Trudinger[J]. Indiana University Mathematics Journal, 1971, 20(11): 1077-1092.	

\bibitem{ref4}Peetre J. Espaces d'interpolation et theorem\'{e} de Soboleff[C]. Annales de l'Institut Fourier, 1966, 16: 279-317.
	
\bibitem{ref3}Pohozaev S I. The Sobolev embedding in the special case $ pl=n $[C]. Proceedings of the technical scientific conference on advances of scientific research 1964-1965, Mathematics Sections, 158-170.	

\bibitem{ref39}Serrin J. Local behavior of solutions of quasi-linear equations[J]. Acta Mathematica, 1964, 111: 248–302.

\bibitem{ref17}Shan W W, Li X M. Compactness of extremals for  Trudinger-Moser functionals on the unit ball in $\mathbb{R}^{2} $[J]. Acta Mathematica Sinica, English Series, 2024, 40(11): 2840-2854.
	
\bibitem{ref8}Struwe M. Critical points of embeddings of $H^{1,n}_0$ into Orlicz spaces[J]. Annales de l'Institut Henri Poincar\'{e} Analyse Non Lin\'{e}aire, 1988, 5(5): 425-464.	
	
\bibitem{ref16}Su X F, Xie R L, Li X M. Compactness of extremals for singular Moser-Trudinger functionals in high dimension[J]. Annali di Matematica Pura ed Applicata (1923-), 2024: 1-23.	

\bibitem{ref51}Talenti G. Elliptic equations and rearrangements[J]. Annali della Scuola Normale Superiore di Pisa-Classe di Scienze, 1976, 3(4): 697-718.

\bibitem{ref40}Tolksdorf P. Regularity for a more general class of qusilinear elliptic equations[J]. Differential Equations, 1984, 51(1): 126-150.

\bibitem{ref5}Trudinger N S. On embeddings into Orlicz space and some applications[J]. Journal of Mathematics and Mechanics, 1967, 17(5): 473-484.

\bibitem{ref36}Wang G F, Xia C. A characterization of the Wulff shape by an overdetermined anisotropic PDE[J].  Archive for Rational Mechanics and Analysis, 2011, 99: 99-115.

\bibitem{ref20}Wang G F, Xia C. Blow-up analysis of a Finsler-Liouville equation in two dimensions[J]. Journal of
Differential Equations, 2012, 252(2): 1668-1700.

\bibitem{ref15}Wang Y M, Yang Y Y. Compactness of extremals for critical singular Trudinger-Moser functions[J]. Journal of Mathematical Analysis and Applications, 2021, 496(2): 124841.

\bibitem{ref26}Xie R L. Anisotropic Moser–Trudinger Inequality Involving $L^{n}$ Norm in the Entire Space $\mathbb{R}^{2}$[J]. Acta Mathematica Sinica, English Series, 2023, 39(12): 2427-2451.	
	
\bibitem{ref44}Xie R L,  Gong H J. A priori estimates and blow-up behavior for solutions of $-Q_{N}u=Ve^{u}$ in 
bounded domain in $\mathbb{R}^{N}$[J]. Science China Mathematics, 2016, 59: 479-492.	
	
\bibitem{ref33}Yang Y Y. Extremal functions for Trudinger-Moser inequalities of Adimurthi-Druet type in dimension two[J]. Journal of Differential Equations, 2015, 258(9): 3161-3193.	

\bibitem{ref13}Yang Y Y, Zhu X B. Blow-up analysis concerning singular Trudinger-Moser inequalities in dimension two[J]. Journal of Functional Analysis, 2017, 272(8): 3347-3374.

\bibitem{ref43}Yang Y Y, Zhu X B. A Trudinger-Moser inequality for conical metric in the unit ball[J]. Archiv der Mathematik, 2019, 112(5): 531-545.

\bibitem{ref2}Yudovich V I. Some estimates connected with integral operators and with solutions of equations[J]. Doklady Akademii Nauk. Russian Academy of Sciences, 1961, 138(4): 805-808.

\bibitem{ref35}Zhou C L, Zhou C Q. Extremal functions of the singular Moser-Trudinger inequality involving the eigenvalue[J]. Journal of Partial Differential Equations, 2018, 31(1): 71-96.

\bibitem{ref22}Zhou C L, Zhou C Q. Moser-Trudinger inequality involving the anisotropic Dirichlet norm $(\int_{\Omega}^{}F^{N}(\nabla u)dx)^{1/N}$ on $W_{0}^{1,N}(\Omega)$[J]. Journal of Functional Analysis, 2019, 276 (9): 2901-2935.

\bibitem{ref25}Zhou C L, Zhou C Q. On the anisotropic Moser-Trudinger inequality for unbounded domains
in $\mathbb{R}^{n}$[J]. Discrete and Continuous Dynamical Systems, 2020, 40(2): 847-881.

\bibitem{ref23}Zhu X B. Remarks on singular Trudinger-Moser type inequality[J]. Communications on Pure and Applied Analysis, 2020, 19(1): 103-112.

\end{thebibliography}

\end{document}